%% file: main_black.tex
\documentclass[review, 11pt]{elsarticle}
\journal{}

\biboptions{authoryear}

\usepackage[english]{babel}
\usepackage{amsmath,amsthm,amsfonts,amssymb,amstext}

\usepackage[dvipsnames]{xcolor} 
\usepackage{mathrsfs} 
\usepackage{natbib}
\usepackage{float} 
\usepackage{hyperref} 
\usepackage{url}
\usepackage{color}
\usepackage{colortbl}

\usepackage{longtable} 
\usepackage{pdflscape} 
\usepackage{booktabs} 
\usepackage{multicol,multirow}

\usepackage[normalem]{ulem}  

\usepackage[margin=2cm]{geometry}

\usepackage{tikz}
\usetikzlibrary{patterns}
\usepackage{tikz-cd} 
\usetikzlibrary{positioning} 
\usetikzlibrary{arrows.meta} 
\usetikzlibrary{arrows.meta}
\tikzset{%
  >={Latex[width=2mm,length=2mm]},
            base/.style = {rectangle, rounded corners, draw=black,
                           minimum width=4cm, minimum height=1cm,
                           text centered, font=\sffamily},
  activityStarts/.style = {base, fill=blue!30},
       startstop/.style = {base, fill=red!30},
    activityRuns/.style = {base, fill=green!30},
         process/.style = {base, minimum width=2.5cm, fill=orange!15,
                           font=\ttfamily},
}

\usepackage{lmodern} 

\usepackage{setspace} 
\usepackage{orcidlink}

\newpageafter{author} 

\newcommand{\printmargins}{
  \noindent
  Left: \the\oddsidemargin \\
  Right: \the\evensidemargin \\
  Top: \the\topmargin \\
  Textwidth: \the\textwidth \\
  Textheight: \the\textheight \\
}

\begin{document}

\setcounter{page}{0} 

\begin{frontmatter}
\title{Optimization-based strategic planning for geographical\\ healthcare accessibility in developing countries: a literature review}

\author[mymainaddress,mysecondaryaddress]{Laura Davila-Pena\corref{mycorrespondingauthor}
\orcidlink{0000-0003-2175-2546}}
\cortext[mycorrespondingauthor]{Corresponding author}
\ead{lauradavila.pena@usc.es}

\author[mysecondaryaddress]{Maria Paola Scaparra\orcidlink{0000-0002-2725-5439}}
\ead{m.p.scaparra@kent.ac.uk}

\author[mythirdaddress]{Dick den Hertog}
\ead{d.denhertog@uva.nl}

\address[mymainaddress]{MODESTYA Research Group, Department of Statistics, Mathematical Analysis and Optimization, Faculty of Mathematics, Universidade de Santiago de Compostela, Campus Vida, 15782 Santiago de Compostela, Spain.}
\address[mysecondaryaddress]{Centre for Logistics and Sustainability Analytics (CeLSA), Department of Analytics, Operations and Systems, Kent Business School, University of Kent, CT2 7PE Canterbury, UK.}
\address[mythirdaddress]{Amsterdam Business School, University of Amsterdam, Amsterdam 1018TV, the Netherlands.}

\begin{abstract} 
Access to healthcare facilities is a critical issue in developing countries, where limited resources and significant challenges hinder progress toward universal health coverage, one of the targets pursued by the United Nations. As a result, the Operations Research (OR) community has become increasingly active in addressing this issue, employing various techniques, particularly optimization. However, making long-term strategic decisions remains challenging, as these often require substantial investments and directly affect large populations. Moreover, much of the existing literature is predominantly theoretical. This paper reviews studies published over the past two decades that use optimization techniques to inform strategic planning decisions aimed at improving health accessibility and specifically contain a case study conducted in a developing country. The review explores the nature of the problems tackled, considering factors such as the level of care, service delivery channels, objectives, criteria, or potential hierarchy, uncertainty, and multi-period settings. Additionally, we discuss the modeling approaches, solution methodologies, and the extent of practical implementation. The goal of this survey is not only to summarize existing work but also to provide a roadmap for future research, offering valuable insights for OR practitioners in the field.

\end{abstract}


\end{frontmatter}


\section{Introduction}\label{sec:intro}

In developing countries, ensuring equitable access to healthcare services remains a major challenge, marked by inadequate infrastructures or geographical, cultural, and economic barriers. The 2030 Agenda for Sustainable Development \citep{UN2015} establishes the 17 Sustainable Development Goals (SDGs), of which several address health-related objectives. Notably, SDG3 concentrates on ``ensuring healthy lives and promoting well-being for all at all ages'' \citep{UN2025}, where target 3.8 strives for universal health coverage (UHC). UHC aims to ensure that all individuals and communities have access to quality health services without incurring financial burdens.

In 2023 the United Nations (UN) reported insufficient progress in expanding UHC \citep{UNDESA2023}, which disproportionately impacts developing nations and affects the achievement of other SDGs \citep{WHO2017}. About 70\% of maternal deaths occur in Sub-Saharan Africa \citep{WHO2019} and Mozambique ranks among the top ten countries worldwide with the highest HIV prevalence \citep{john2020predictors}. In Mali, malaria constitutes a primary cause of morbidity and mortality among children under the age of five \citep{oliphant2022improving} and mean rates of treatment coverage are under 50\% in countries like Bangladesh \citep{rahman_progress_2017} or Dominican Republic \citep{mathon2018cross-border}. Naturally, disparities are further exacerbated among diverse socioeconomic groups. For example, in the rural regions of Northern Brazil, life expectancy is notably lower compared to more populated areas, and substantial differences exist in the under-five mortality rate \citep{benevenuto2019assessing}. 
%

There are numerous factors that drastically impact the health of the most vulnerable populations, resulting in consequences as described above. Among these causes are the distance to healthcare facilities (\citealp{feikin_impact_2009}; \citealp{mcguire2021effect}; \citealp{Lefevre2025}), the scarcity of transportation means, particularly in rural or mountainous regions \citep{agbenyo2017accessibility}, poor road conditions (\citealp{agbenyo2017accessibility}; \citealp{Peprah2020does}), the absence of an organized and centralized ambulance system (\citealp{boutilier2020ambulance}; \citealp{nadar2023adaptive}), or the lack of a sufficient referral process between different-level medical services \citep{khodaparasti2017enhancing}. Additionally, the implementation of inadequate vaccination procedures contributes to the rise of the non-immunized population, fostering the spread of diseases and epidemics \citep{chen2014planning}. Furthermore, disruptive events such as floods (\citealp{cook2019integrating}; \citealp{de2022vaccine}; \citealp{van2022optimizing}), geological hazards (\citealp{esposito_amideo_optimising_2019}; \citealp{yu2023comprehensive}), and fires \citep{rathore2022sustainable} cause systems in these countries to collapse, with inevitable consequences for access to healthcare infrastructure. On a sociological level, distrust in the healthcare system \citep{onwujekwe_corruption_2019}, lack of education and information \citep{thongmixay2019perceived}, as well as unfavorable laws enacted by certain governments \citep{ferguson2022}, also play a considerable role. Lastly, the significant population growth and the increasing aging society (\citealp{khodaparasti2018multi-period};  \citealp{dogan_model_2020}) 
add additional pressure to healthcare systems and exacerbate challenges in delivering health services. 

The scientific community is increasingly committed to tackling challenges in developing countries and advancing the SDGs. \cite{orhan2022analytics} reviewed research papers on business analytics in these countries, categorizing them by methodology, application area, and SDG relevance. Their survey found optimization to be the most used approach for decision-making problems. However, among the optimization-based studies they identified, only a small subset focused on healthcare (\citealp{murawski2009improving}; \citealp{smith2009planning}; \citealp{cocking2012improving}; \citealp{andrade2015abc}; \citealp{chaiwuttisak2016location}; \citealp{smith2018siting}; \citealp{buyuktahtakin2018new}; \citealp{savacser2019organ}; \citealp{yuniaristanto2019assessing}) and even fewer---the first four---addressed strategic planning to improve healthcare accessibility. 
Our work builds on this observation by providing a focused review of optimization-based models that support long-term strategic decisions aimed at enhancing accessibility to health services in developing countries. While developed countries often face challenges related to efficiency or equity in already well-functioning systems, operational challenges in developing countries are typically more severe and multifaceted, as showcased above. In such settings, optimization techniques are especially suitable because they rigorously support decision-making under objective functions and constraints, systematically evaluate trade-offs, and provide actionable guidance for planning in environments where ad-hoc or manual methods are often the norm.

Combining all the above, this survey reviews and classifies research articles from the last two decades that propose strategic planning decisions based on optimization, enhance patients' accessibility to healthcare services, and include an application in developing countries tailored to their specific features. Section~\ref{sec:scope} further elaborates on the scope of the survey. The goal is to identify prevalent optimization techniques and potential gaps in the literature. Specifically, 56 papers meet these criteria. We analyze key problem features, proposed modeling and solution approaches, and the extent to which the findings have been implemented in practice. Notably, only four papers report real-world implementation. Given the vast potential of optimization to solve problems in this area, we present a research roadmap, highlighting essential practical features, novel optimization methods, and tools for successful implementation and measurable impact.

To the best of our knowledge, this is the first survey with these characteristics, specifically applied in developing countries and addressing various strategic actions aimed at improving healthcare access for their populations. Previous surveys focusing on developing countries include \cite{rahman_use_2000}, which explores location and allocation models in healthcare, and \cite{white2011developing}, which revises the broader literature on OR techniques but does not specifically focus on the health domain. In contrast, general reviews of OR applications in healthcare, which might contain overlapping literature to this survey, include \cite{brailsford2011}, \cite{rais2011operations}, and  \cite{belien2024}, although they do not focus on developing countries. Additionally, \cite{ahmadi2017survey} and \cite{gunecs2019location} provide surveys on healthcare facility location up to their respective dates of publication.

The organization of this paper is as follows. Section~\ref{sec:scope} presents the scope of the review. Section~\ref{sec:methodology} describes the adopted methodology for the literature search. Section~\ref{sec:main} analyzes the selected papers according to the problem features, the modeling and solution approaches, and the practical implementation of the proposed models. Section~\ref{sec:roadmap} offers a roadmap for future research that follows the structure of the previous section. The paper finishes with some concluding remarks in Section~\ref{sec:conclusions}.

\section{Scope of the review}\label{sec:scope}

This review analyzes optimization-based strategic interventions aimed at improving the accessibility to health services in developing countries. The term \textit{accessibility} encompasses a range of interpretations, spanning broader aspects such as the influence of education availability on healthcare-seeking behaviors to more straightforward considerations, such as spatial proximity for physical access to medical attention \citep{agbenyo2017accessibility}. Our study narrows its focus to the fundamental concept of accessibility: geographical accessibility and the ease with which end-users can access healthcare services or facilities (e.g., in terms of distance, time, or available transport means/network) to receive assistance. 
Thus, we concentrate on interventions that directly impact the ability of individuals to access medical care. By \textit{end-users}, we refer to patients or population groups in need of healthcare. 

The review is restricted to strategic planning decisions that involve a real-world case study in a \textit{developing country}, defined here using the DAC list of ODA Recipients \citep{ODA2022}, which includes least developed countries (LDCs) and low- and middle-income countries (LMICs) according to gross national incomes per capita from the World Bank. 
By \textit{strategic planning decisions}, we refer to long-term, high-level decisions such as where to locate health facilities, how to design or expand service delivery systems, and which infrastructure developments should be carried out to improve physical accessibility to health services (e.g., road improvements leading to health facilities). These decisions typically set the foundation for subsequent tactical (e.g., health products or capacity allocation, staffing levels) and operational (e.g., daily routing, workforce scheduling) actions. While strategic planning decisions can arise in both stable and crisis contexts, this review focuses on those made under routine (non-emergency) circumstances. 
Thus, we deliberately exclude papers dealing with disaster response or pandemic preparedness, even when they involve strategic decisions. These topics have been extensively dealt with in the literature. For instance, \cite{galindo2013review} and \cite{ozdamar2015models} survey disaster operations management and humanitarian logistics studies, including long-term planning phases that could involve health systems, while \cite{boonmee2017facility} reviews facility location models in emergency humanitarian logistics. Similarly, pandemic preparedness models, including strategic planning during COVID-19, have been recently covered in \cite{dey2024optimization}.

Finally, while both healthcare services and health product supply chains influence accessibility, this review concentrates on \textit{health services}. 
Papers focusing primarily on the supply chain of health products are excluded, unless they directly address patient access to healthcare facilities or services. For example, studies on vaccine supply chains (see \citealp{duijzer2018literature}, and \citealp{deboeck2020vaccine}, for comprehensive reviews) may be included when they consider the final stage of delivery involving patient immunization, thereby constituting a health service.

In summary, to select the papers for inclusion in our review, the following criteria were established: i) The paper addresses the improvement of geographical accessibility to healthcare services; ii) The paper includes a case study conducted in a developing country, as defined above; iii) The paper uses optimization-based  methods; iv) The paper focuses on long-term strategic planning decisions in non-emergency contexts; v) The paper targets access to health services by patients or population groups.

\section{Search methodology}\label{sec:methodology}

The methodology employed in this review involved systematic search and selection procedures to identify relevant literature within the domain of healthcare accessibility improvement through optimization, with a focus on developing countries. The search process encompassed the period from January 2003 to June 2025 and was conducted across INFORMS journals, ScienceDirect, and Scopus databases.

The search strategy involved constructing Boolean strings tailored to each database's capabilities and incorporating terms such as ``accessibility'', ``Operations Research'', ``optimization'', and ``developing countries'', along with several variations or combinations. We intentionally did not include keywords related to ``strategic'' in our search strings, as many long-term decisions (such as the location of new health centers) are strategic in nature, even if not explicitly labeled as such in the original work. Instead, we applied our own criteria, together with how the original studies framed their decisions when available, to assess whether a paper addressed a long-term strategic planning problem. Additional details on the search strategy can be found in the Supplementary Material.

Figure~\ref{fig:PRISMA} presents a PRISMA flow diagram \citep{PRISMA2021} summarizing the literature search process. More than 2,500 records were identified and initially screened based on their titles and keywords to determine potential relevance. Over 200 records were sought for retrieval, for which abstracts and conclusions were reviewed. Approximately half of these underwent full-text examination to assess their eligibility. To select the papers for inclusion in our review, the criteria established in Section~\ref{sec:scope} were considered. Reasons for exclusion included: a focus on tactical or operational decisions rather than strategic planning; the absence of a real-world case study, or a case study conducted in a developed country; the use of purely descriptive methodologies (e.g., geographic information systems); a focus on health products rather than health services; or strategic decisions that did not directly impact patients (e.g., location of intermediate storage hubs for vaccines rather than vaccination delivery points). 
Finally, 56 papers were selected for comprehensive analysis and review. For a brief bibliometric analysis of the selected papers, we refer the reader to the Supplementary Material.

\begin{figure}[ht!]
    \centering
    \includegraphics[width=\textwidth]{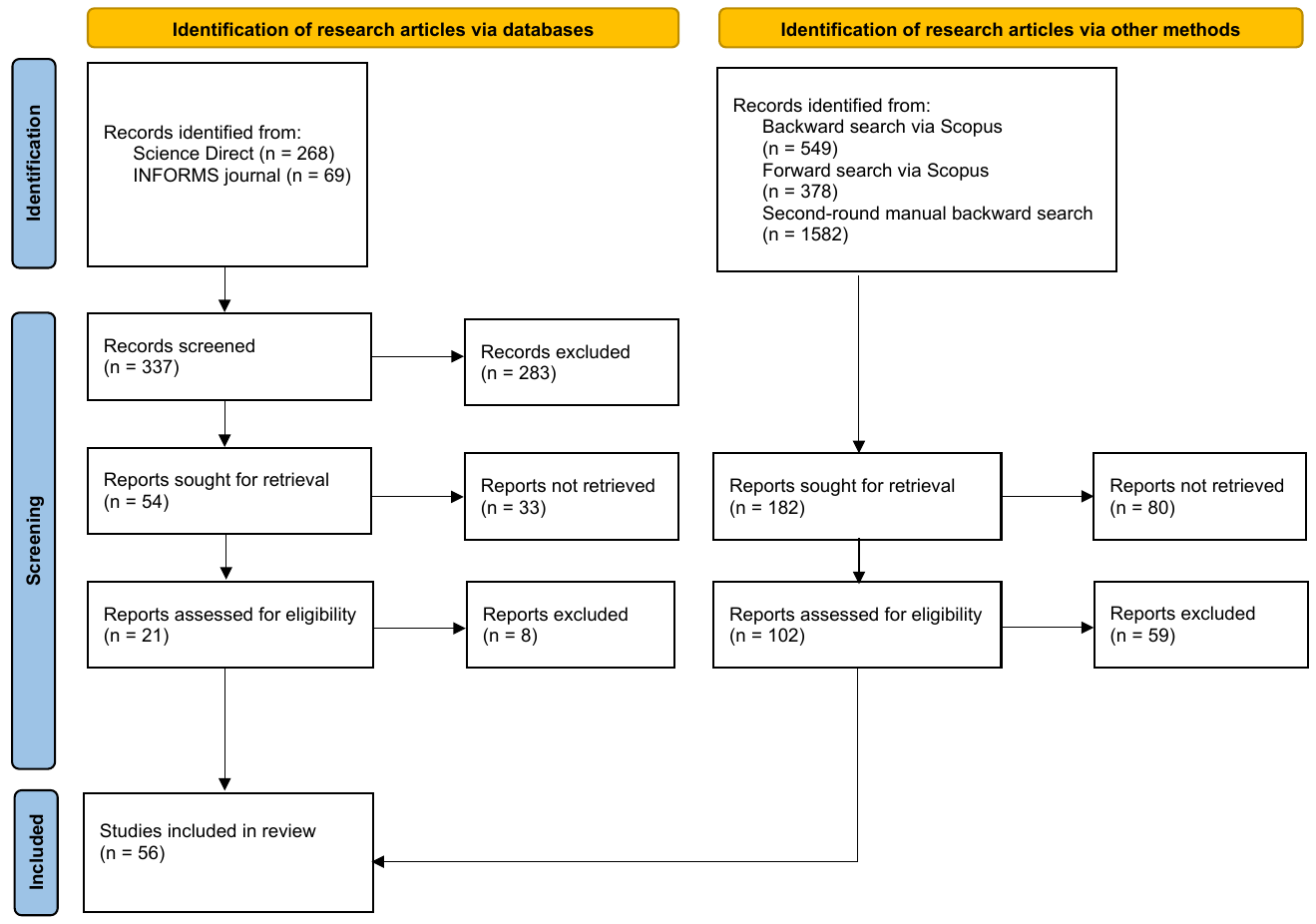}
    \caption{PRISMA flow diagram for this review paper.}
    \label{fig:PRISMA}
\end{figure}

\section{Selected papers analysis}\label{sec:main}
In the subsequent subsections, we present a comprehensive analysis of the literature focusing on the features of the addressed problems, the modeling and solution approaches adopted, and the specific characteristics of case studies regarding their practical implementation. Within each of these subsections, several pertinent aspects are discussed, accompanied by tables aimed at offering additional insights.

\subsection{Problem features}\label{subsec:problem_features}
We examine the level of healthcare under consideration, the delivery channel for such care reaching patients, the pursued objectives and high-level criteria considered, the potential presence of hierarchy, uncertainty, or multi-period settings, and the geographical area and country of study. All these aspects are summarized in Table~\ref{tab:problem_features}, where references are listed in chronological order.

\subsubsection{Level of care}\label{subsec:type_of_care}
The healthcare systems of each country heavily depend on its specific characteristics, and generally, the care offered by each medical service tends to vary from one country to another. Thus, categorizing the level of care accurately is not a straightforward task, and at times the line between one and another is blurred. In an attempt to faithfully depict the level of care that each of the selected articles has focused on, we present the third column of Table~\ref{tab:problem_features}.



\begin{table}[ht!]\footnotesize
\centering
\caption{Problem features for the selected papers.}\label{tab:problem_features}
\resizebox{1.05\linewidth}{!}{
\hspace{-1.2cm}\begin{tabular}{p{0.5cm}p{5.9cm} llll lll lll}
\noalign{\smallskip}\hline\noalign{\smallskip} 
& & \multirow{2}{*}{\parbox{1.5cm}{Level of care}} & \multirow{2}{*}{\parbox{1.5cm}{Service delivery channel}} & \multicolumn{2}{l}{\multirow{2}{*}{\parbox{1.6cm}{Types of objectives}}} &  Criteria &  Hierarchy & Uncertainty & \multirow{2}{*}{\parbox{1.2cm}{Multi-period}} & \multirow{2}{*}{\parbox{1.8cm}{Geographical area}} & Country\\
No. & Reference & & & & & & & & & \\
   \cmidrule(l){5-6} 
&  & & & Patient & System & & & & &  \\
\noalign{\smallskip}\hline\noalign{\smallskip}
    1 &  \cite{galvao2006load}  &  MHC   &   FF &  D  &  Im 
    &  Ac, Eq &   \checkmark &  &     &  Municipality (U)     &  Brazil \\[0.2em]
    \arrayrulecolor{gray!30}\hline

    2 & \cite{doerner2007multicriteria}  & PHC   &  MF &  Cov, D  & C &  Ac, Eff, Eq &   &  & & Region (R, U) 
    & Senegal \\[0.2em]
     \arrayrulecolor{gray!30}\hline
    
    3 & \cite{murawski2009improving}     & PHC   &  FF  &  Cov & & Ac &   &   &  & District (R)  & Ghana\\[0.2em]
     \arrayrulecolor{gray!30}\hline
     
    4 & \cite{smith2009planning}     & CHC  & IHW & Cov, D & nW & Eff, Eq & \checkmark &  &  & District (R) & India\\[0.2em] \arrayrulecolor{gray!30}\hline
    
    5 & \cite{cocking2012improving}       & PHC  & FF &  C & & Ac  &   &  & & Health district (R) & Burkina Faso\\[0.2em] \arrayrulecolor{gray!30}\hline
    
    6 & \cite{shariff2012location}       & PHC &  FF &  Cov &   & Ac, Eff  &   &  & & \multirow{2}{*}{\parbox{3cm}{District (R, U) \\Sub-district (R)}} & Malaysia \\
    & & & & & & & & & & \\[0.2em] \arrayrulecolor{gray!30}\hline
    
    7 & \cite{ghaderi2013modeling}     &  PHC & FF & C & C  & Ac  &   &  & \checkmark & Province (R, U) 
    & Iran\\[0.2em] \arrayrulecolor{gray!30}\hline

    8 &  \cite{smith2013bicriteria}  &  G   &   FF &  Cov, D  &   &  Eff, Eq &   \checkmark &  & &  Region (R)   &  India \\[0.2em]
    \arrayrulecolor{gray!30}\hline

    9 & \cite{gunecs2014}     & PHC & FF   & Ac, Cov, D 
    &  & Ac, Eq  &   &  &  & Province (R, U) & Turkey \\[0.2em] \arrayrulecolor{gray!30}\hline

    10 & \cite{kunkel2014optimal}   & CHC & IHW & C, D & C  & Ac, Eff, Eq  &   &  &  & \multirow{2}{*}{\parbox{3cm}{District (R) \\Country-wide}} & Malawi \\
    & & & & & & & & & & \\[0.2em] \arrayrulecolor{gray!30}\hline

    11 & \cite{mohammadi2014design}      & G
    & FF & C, T  & C   &  Ac, Eff
    &  & \checkmark & & Province (U) & Iran\\[0.2em] \arrayrulecolor{gray!30}\hline

    12 & \cite{yao2014locational}      &  SRH  & FF &  D & & Ac  &   &  & & Districts (R) & Mozambique\\[0.2em] \arrayrulecolor{gray!30}\hline

    13 & \cite{zahiri2014multi-objective}      & OT & FF   &  & C, T & Ac, Eff  &   &   \checkmark & \checkmark & Country-wide & Iran\\[0.2em] \arrayrulecolor{gray!30}\hline

    14  & \cite{zahiri2014robust}      & OT & FF  &  & C  & Eff   &   &  \checkmark & \checkmark & Province (U) & Iran\\[0.2em] \arrayrulecolor{gray!30}\hline

    15  & \cite{andrade2015abc}     &  EMS   & MF & Cov & T & Eff  &   &  \checkmark & \checkmark & City (U) & Brazil\\[0.2em] \arrayrulecolor{gray!30}\hline
    
    16  & \cite{beheshtifar_multiobjective_2015}    & PHC & FF & C, D & C, D  & Ac, Eq, S-s    &   &  & & Region (U)  & Iran\\[0.2em] \arrayrulecolor{gray!30}\hline

    17  & \cite{shishebori2015robust}      &  PHC & FF &  & C  & Ac, Eff, Eq  &   &  \checkmark & & Towns 
    (U) & Iran\\[0.2em] \arrayrulecolor{gray!30}\hline
    
    18  & \cite{lim2016coverage}    &  V & MF & Cov &  & Ac &   &  & & State (R) & India\\[0.2em] \arrayrulecolor{gray!30}\hline
    
    19  & \cite{nunez2016column}    & PHC+PrevHC  &  MF &  Cov, d, Eq &   & Ac, Eq  &   &  & & Transport corridors  & SSA\\[0.2em] \arrayrulecolor{gray!30}\hline
    
    20  & \cite{vonachen2016optimizing}    & CHC & IHW  & & C, nW  &  Ac &   &   & & Districts (R) & Liberia\\[0.2em] \arrayrulecolor{gray!30}\hline
    
   21  & \cite{khodaparasti2017enhancing}  & CHC &  FF & T & D, Eff 
   &  Ac, Eff, Eq &  \checkmark &  & & City (U) & Iran\\[0.2em] \arrayrulecolor{gray!30}\hline

    22 & \cite{zarrinpoor2017design}  &  G   &   FF &   C &  C &  Ac, Eff &   \checkmark &   \checkmark & &  Province (R,U)     &  Iran \\[0.2em]
    \arrayrulecolor{gray!30}\hline

    23  & \cite{khodaparasti2018multi-period}    & ECS & FF    & Cov  &   & Ac &   & \checkmark  & \checkmark & City (U) & Iran\\[0.2em] \arrayrulecolor{gray!30}\hline
    
    24  & \cite{mousazadeh_accessible_2018}    & G &  FF &  Cov & In, SL  & Eq, Sta  &  \checkmark & \checkmark & \checkmark & City (U) & Iran\\[0.2em] \arrayrulecolor{gray!30}\hline

    25 &  \cite{mousazadeh2018health}  &   G   &   FF &  T  &  C &  Ac, Eff &  \checkmark &   \checkmark &  & City (U)     &  Iran \\[0.2em]
    \arrayrulecolor{gray!30}\hline

    26  & \cite{wang2018healthcare}    &  G & FF & D, EB & C, Q  &  Ac, Eff &  \checkmark & \checkmark & & County (R) & China\\[0.2em] \arrayrulecolor{gray!30}\hline
    
    27  & \cite{zarrinpoor2018design}    &  G  &  FF & C  &  C & Ac, Eff &  \checkmark & \checkmark  & & Province (R, U) & Iran\\[0.2em] \arrayrulecolor{gray!30}\hline

    28  & \cite{cherkesly2019community}   & CHC &  IHW &   & C & Ac, Eff  &   &   & & Districts (R) & Liberia\\[0.2em] \arrayrulecolor{gray!30}\hline
    
    29  & \cite{boutilier2020ambulance}   & EMS  &  MF &  & T & Ac, Eff, Eq  &   & \checkmark  & & City (U) & Bangladesh\\[0.2em] \arrayrulecolor{gray!30}\hline

    30  & \cite{de2020roadside}    &  PHC+PrevHC & MF & d 
    &  Effect &  Ac, Eff, Eq &   &  & & Transport corridors & SSA\\[0.2em] \arrayrulecolor{gray!30}\hline

    31  & \cite{dogan_model_2020}    &  PrevHC & FF  & Cov & C, T  & Ac, Eff  &   & \checkmark & \checkmark & City (U) & Turkey\\[0.2em] \arrayrulecolor{gray!30}\hline
    
    32  & \cite{motallebi2020short-term}    & G & FF  &  T & C & Ac, Eff & \checkmark &  \checkmark  & & Province (U) & Iran\\[0.2em] \arrayrulecolor{gray!30}\hline

    33  & \cite{taymaz2020healthcare}     &  PHC+PrevHC & MF & Cov  &   & Ac  &   &  \checkmark & & Transport corridors & SEA\\[0.2em] \arrayrulecolor{gray!30}\hline

    34  & \cite{almeida2021two-step}    &  SHC  &  FF &  d, D & D  & Ac, Eff, Eq  &   &  & & State (R, U) & Brazil\\[0.2em] \arrayrulecolor{gray!30}\hline
    
    35  & \cite{flores2021optimizing}    &  PHC  &  FF & Cov  & & Ac, Eff, Eq  &   &   & & City (R, U) & Philippines\\[0.2em] \arrayrulecolor{gray!30}\hline

    36  & \cite{srivastava2021strengthening}  & V &  FF  & Cov &   & Eq   &   &   & &  Districts (R, U) & India\\[0.2em] \arrayrulecolor{gray!30}\hline
        
    37 & \cite{yang2021outreach}    & V  & MF &  & C & Ac &  &  \checkmark & \checkmark & Regions (R) & SSA\\[0.2em] \arrayrulecolor{gray!30}\hline

    38 &  \cite{alban2022resource}  &  SRH   &   MF &  d  &   &  Ac, Effect &   &   \checkmark & &  Country-wide     &  Uganda \\[0.2em]
    \arrayrulecolor{gray!30}\hline

    39 & \cite{chen2022multi-objective}   &  EMS & FF & Cov, T  &  C &  Ac, Eff &  &  & & County (R) 
    & China\\[0.2em] \arrayrulecolor{gray!30}\hline

    40 & \cite{chouksey2022hierarchical}   & MHC &  FF & C  & C  & Ac, Av & \checkmark &   & & Administrative unit (R) & 
    India\\[0.2em] \arrayrulecolor{gray!30}\hline

    41  & \cite{chouksey2022optimisation}    & MHC &  FF &  C & C  & Ac, Av
    &  \checkmark & \checkmark & & District (R, U) & India\\[0.2em] \arrayrulecolor{gray!30}\hline

   	42  & \cite{de2022optimization}    & MRI  &  FF & d, D  & C  & Ac, Eff, Eq  &   &  & & Country-wide & Brazil\\[0.2em] \arrayrulecolor{gray!30}\hline

    43  & \cite{elorza2022assessing}    & PHC  & FF & Cov &   &  Ac, Eq & &  & & City (U) & Argentina\\[0.2em] \arrayrulecolor{gray!30}\hline

    44 &  \cite{karakaya2022biobjective}  &  MHC   &   FF &  C  &  C &  Ac, Eff &   \checkmark &    & &  Region (R, U)     & Turkey \\[0.2em]
    \arrayrulecolor{gray!30}\hline

    45  & \cite{mendoza-gomez2022location}   & PHC  &  FF   & Cov, D  & C  & Ac &  &  & & Region (R, U)
    & Mexico\\[0.2em] \arrayrulecolor{gray!30}\hline

    46  & \cite{mendoza-gomez2022regionalization}    & PHC &  FF &  D & & Ac    &   &   & & State (R, U) 
    &   Mexico \\[0.2em] \arrayrulecolor{gray!30}\hline

    47  & \cite{pourrezaie-khaligh2022fix-and-optimize} & PHC
    & FF  &  Ac, C & C, En  &  Ac, Eq &  & & \checkmark &  Province (R, U)
    &   Iran\\[0.2em] \arrayrulecolor{gray!30}\hline
    
   	48  & \cite{rouhani2022robust}     &  OT &   FF  &  C, D, T & C, T  & Ac, Eq  & \checkmark & \checkmark & \checkmark &  Provinces (U) & Iran\\[0.2em] \arrayrulecolor{gray!30}\hline

    49  & \cite{kumar2023location}     & PHC &  FF &  Cov, d & & Ac, Eq  &   &  & & District (U) & India\\[0.2em] \arrayrulecolor{gray!30}\hline

    50 &  \cite{decampos2024multi}  &  MS   &   FF &  Cov, D, Eq  &   &  Ac, Eq &    &    & &  States (U)     &  Brazil \\[0.2em]
    \arrayrulecolor{gray!30}\hline

    51 &  \cite{defreitas2024integrated}  &  NS   &   FF &   &  C &  Ac, Eff &   &    & &  Country-wide     &  Brazil \\[0.2em]
    \arrayrulecolor{gray!30}\hline

    52  & \cite{eksioglu2024designing}   &  V & FF &  Cov & nF  & Ac  &  &  & \checkmark &  Regions (R) & Niger\\[0.2em] \arrayrulecolor{gray!30}\hline

    53  & \cite{heyns2024optimisation}    &  G
    &  FF & Ac  & C  & Ac & \checkmark & & & Province (R)    & Nepal\\
    \arrayrulecolor{gray!30}\hline
 
    54 &  \cite{mendoza2024maximal}  &  SHC   &   FF &  Cov  &  &  Ac &   &   & &  Country-wide     &  Mexico \\[0.2em]
    \arrayrulecolor{gray!30}\hline

    55 &  \cite{mendoza2024regionalization}  &  PHC   &  FF &  D  &   &  Ac &    &   & &  State (R, U)     &  Mexico \\[0.2em]
    \arrayrulecolor{gray!30}\hline
  
    56 &  \cite{chouksey2024accelerated}  &  MHC   &   FF &  C  &  C & Ac, Av &   \checkmark &    & \checkmark & District (R)    &  India \\

     \arrayrulecolor{black}
\noalign{\smallskip}\hline\noalign{\smallskip}
    \multicolumn{12}{l}{\parbox{1.55\linewidth}{\textbf{Level of care}: MHC (maternal healthcare), PHC (primary healthcare), CHC (community healthcare), G (general), SRH (sexual and reproductive health), OT (organ transplant), EMS (emergency medical services), V (vaccination), PrevHC (preventive healthcare), ECS (elderly care services), SHC (secondary healthcare), MRI (magnetic resonance imaging), MS (mammography screening), NS (newborn screening);
    \textbf{Service delivery channel}: FF (fixed facility), MF (mobile facility), IHW (itinerant health worker); 
    \textbf{Types of objectives}: D (distance), Im (Imbalanced facility load), Cov (coverage), C (cost), nW (number of workers), Ac (accessibility), T (time), d (served demand), Eq (equity), Eff (efficiency), In (instability), SL (service level), EB (environmental benefits), Q (quality), Effect (effectiveness), En (envy), nF (number of facilities); 
    \textbf{Criteria}: Ac (accessibility), Eq (equity), Eff (efficiency), S-s (site-suitability), Sta (Stability), Effect (effectiveness), Av (availability);  
    \textbf{Geographical area}: U (urban), R (rural); 
    \textbf{Country}: for those categorized as SSA (Sub-Saharan Africa) or SEA (South-East Africa), the specific countries are not explicitly mentioned.}}
    \end{tabular}}
\end{table}


The literature is heavily dominated by studies focusing on primary care, reflecting its central role as the first point of contact and its importance in prevention and early intervention, particularly in resource-constrained settings. Several of the analyzed papers (15 out of 56) concentrate on primary healthcare (PHC) in general, aiming to improve spatial accessibility to basic healthcare facilities, often in rural or underserved areas across Africa and Asia (\citealp{doerner2007multicriteria, murawski2009improving, cocking2012improving, shariff2012location, ghaderi2013modeling, gunecs2014, beheshtifar_multiobjective_2015, shishebori2015robust, flores2021optimizing, pourrezaie-khaligh2022fix-and-optimize, kumar2023location}). Some works focus on locating, upgrading, or redesigning primary care facilities within South American healthcare systems, either by introducing complementary services such as nutritional advice, dental care,
mental health services, clinical analysis, or imaging services  (\citealp{mendoza-gomez2022location,mendoza-gomez2022regionalization,mendoza2024regionalization}) or by reallocating basic services within existing primary care centers (\citealp{elorza2022assessing}).

Beyond PHC in general, other studies (16 out of 56) focus on specific primary-level services. 
For instance, some address preventive healthcare (PrevHC), including the location of early diagnosis and cancer screening centers (\citealp{dogan_model_2020}). Other works tackle the location of roadside healthcare facilities along major transport corridors in Africa (\citealp{nunez2016column, de2020roadside, taymaz2020healthcare}) that provide truck drivers with basic care and prevention services. Some papers address the installation of newborn screening (NS) facilities and equipment (\citealp{defreitas2024integrated}), or sexual and reproductive health (SRH) and related services, modeling the relocation and expansion of HIV service provision (\citealp{yao2014locational}) or the allocation of mobile family-planning healthcare units to underserved areas (\citealp{alban2022resource}).
A further stream concentrates on vaccination (V), aiming to improve immunization rates among the population by appropriately locating outreach posts or session sites to deliver routine vaccines (\citealp{lim2016coverage,srivastava2021strengthening,yang2021outreach,eksioglu2024designing}). Five contributions model community healthcare (CHC) programs, focusing on the deployment of community health workers, surveillance assistants, or community-based organizations to increase individuals' access to medical care and foster better participation and greater involvement  (\citealp{smith2009planning,kunkel2014optimal,vonachen2016optimizing,cherkesly2019community,khodaparasti2017enhancing}).

All of the above highlights the researchers' interest in addressing primary care issues in developing countries. Fewer works (13 out of 56) focus on secondary healthcare (SHC) broadly (\citealp{almeida2021two-step, mendoza2024maximal}) or specialized care, including magnetic resonance imaging (MRI) (\citealp{de2022optimization}), mammography screening (MS) (\citealp{decampos2024multi}), emergency medical services (EMS) (\citealp{andrade2015abc, boutilier2020ambulance, chen2022multi-objective}), elderly care services (ECS) (\citealp{khodaparasti2018multi-period}), or maternal healthcare (MHC) (\citealp{galvao2006load, chouksey2022hierarchical, chouksey2022optimisation, chouksey2024accelerated, karakaya2022biobjective}). The latter, more represented, group seeks to determine the optimal number and location of health facilities and the assignment of expectant mothers to these centers, where services range from basic maternal care to assistance with childbirth, including unplanned Cesarean sections or neonatal care. Only a small number of works (3 out of 56) consider advanced specialized treatment, specifically addressing the improvement of the efficiency of the organ transplant (OT) network in Iran, where perishability, logistics, and coordination among centers and transportation modes are taken into account (\citealp{zahiri2014multi-objective,zahiri2014robust,rouhani2022robust}).

The remaining selected articles (9 out of 56) either analyze accessibility to all types of care or do not specify the services offered at the healthcare facilities under study. These are classified as general (G) studies in Table~\ref{tab:problem_features} (\citealp{smith2013bicriteria, mohammadi2014design, zarrinpoor2017design, zarrinpoor2018design, mousazadeh_accessible_2018, mousazadeh2018health, wang2018healthcare, motallebi2020short-term,heyns2024optimisation}).

Overall, the literature's strong focus on improving access to primary healthcare underscores its role as the cornerstone of health systems in developing countries. The comparatively limited attention to secondary and specialized care highlights opportunities for future research, particularly in modeling interactions across care levels or facilities, referral pathways, and the integration of primary and specialized services within unified, multi-service planning frameworks.

\subsubsection{Service delivery channel}\label{subsec:service_delivery_channel}
Access to healthcare services can be provided through different delivery channels, depending on whether patients travel to healthcare facilities or medical services are brought closer to the population.  
As depicted in the fourth column of Table~\ref{tab:problem_features}, most papers (43 out of 56) use fixed facilities (FF) to deliver healthcare services, whereas fewer studies investigate mobile facilities (MF) (9 papers) or itinerant health workers (IHW) (4 papers). 

The predominance of fixed facilities reflects the traditional organization of healthcare delivery, where individuals receive healthcare by attending a medical center in person, whether it be a clinic (e.g., \citealp{murawski2009improving}; \citealp{beheshtifar_multiobjective_2015}; \citealp{kumar2023location}), a hospital (e.g., \citealp{zahiri2014multi-objective}; \citealp{zarrinpoor2018design}; \citealp{decampos2024multi}), a transplant center (e.g., \citealp{zahiri2014multi-objective,zahiri2014robust}), a nursing home \citep{khodaparasti2018multi-period}, a vaccination center (e.g., \citealp{srivastava2021strengthening}; \citealp{eksioglu2024designing}), or an emergency station \citep{chen2022multi-objective}, among other instances. This prevalence highlights that most optimization-based planning efforts continue to assume stationary infrastructure as the backbone of healthcare provision.

A smaller but relevant stream of 9 studies departs from this paradigm by explicitly modeling healthcare delivery through mobile facilities. In the domain of emergency services, the use of medicalized vehicles to assist the population is a common practice, jointly investigating the optimal placement of ambulance outposts, alongside the allocation, repositioning or routing of these vehicles (\citealp{andrade2015abc,boutilier2020ambulance}). In contexts characterized by high budgetary constraints, some governments have endeavored to enhance their inhabitants' access to healthcare by implementing mobile units that visit hard-to-reach locations to provide basic healthcare or family planning services (e.g., \citealp{doerner2007multicriteria, alban2022resource}). A related but distinct approach establishes walk-in clinics along major trucking routes for the mobile population in various African countries (\citealp{nunez2016column, de2020roadside, taymaz2020healthcare}). Thus, although these roadside centers cannot strictly be classified as mobile units (given the absence of routing decisions), the location problem addressed accounts for the dynamic nature of truck drivers or sex workers who might need medical care. Meanwhile, some works (\citealp{lim2016coverage,yang2021outreach}) adopt the outreach scheme, where a medical team transports vaccines to remote locations and establishes mobile clinics for immunization.

As mentioned in the previous section, community health programs play a pivotal role in bringing healthcare closer to the most marginalized people in these countries, often involving the mobilization of health workers to underserved areas. This is the case in \cite{smith2009planning}, \cite{kunkel2014optimal}, \cite{vonachen2016optimizing}, and \cite{cherkesly2019community}, which contain considerations regarding the location or routing of healthcare personnel in diverse rural settings of India, Malawi, and Liberia.

Overall, the dominance of fixed-facility solutions and the relatively low consideration of mobile or hybrid approaches highlight opportunities for future research, particularly in designing more flexible and adaptive service delivery systems for rural or infrastructure-poor areas.

\subsubsection{Types of objectives and criteria}\label{subsec:obj_crit}
In the pursuit of enhancing healthcare access in developing countries, strategic planning plays a pivotal role in delineating objectives and guiding decision-making processes. Understanding the diversity of objectives involved within this context is critical for devising effective approaches to address the multifaceted challenges prevalent in healthcare delivery systems.

As can be seen in columns 5 and 6 of Table~\ref{tab:problem_features}, we categorize the objectives identified in the selected literature into two overarching groups: patient-oriented and system-oriented objectives. Patient-oriented objectives focus on improving outcomes from the perspective of the individual seeking healthcare services, while system-oriented objectives entail broader systemic enhancements aimed at optimizing healthcare delivery.
Within these classifications, most of the objectives are further delineated based on underlying high-level criteria, presented in column 7 of Table~\ref{tab:problem_features}, which encompass efficiency (Eff), accessibility (Ac), or equity (Eq), among others. Generally, efficiency aspects revolve around optimizing resource utilization and minimizing costs, while accessibility objectives prioritize reducing barriers such as distance or travel time. Equity criteria, on the other hand, seek to ensure fair and equitable distribution of healthcare services across diverse population segments. 

The problem studied by \cite{doerner2007multicriteria} is a clear example in which the criteria mentioned above are modeled as three objectives. From the user's point of view, the authors propose to minimize the average distance a patient has to walk to reach the nearest stop of the mobile medical facility, as well as to maximize coverage, measured as the percentage of the population living within a predefined maximum walking distance to the medical delivery point. They relate these objectives to accessibility and equity, respectively. On the system side, they aim to minimize staff ineffectiveness, measured as the ratio between non-productive medical time and service time, thus seeking to optimize the time spent providing medical treatment. This is linked to improving the economic efficiency of the system.

Both the different objectives and the criteria tackled by the selected papers may have different definitions and interpretations. From the patient perspective, accessibility- and equity-driven objectives dominate the literature, with coverage-based (Cov) formulations being particularly prevalent (22 out of 56). Most studies define coverage in terms of spatial proximity, aiming to maximize the number of people within a certain distance or travel time from the location where care is provided (\citealp{murawski2009improving, shariff2012location, gunecs2014, lim2016coverage, flores2021optimizing, srivastava2021strengthening, elorza2022assessing, mendoza-gomez2022location, kumar2023location, decampos2024multi, eksioglu2024designing}), or, similarly, to minimize the number of uncovered patients (\citealp{khodaparasti2018multi-period, chen2022multi-objective}). 
Beyond these standard formulations, several studies adopt alternative interpretations of coverage that emphasize balance allocations, continuity, or intensity of access. For instance, \cite{smith2009planning, smith2013bicriteria} minimize deviations from predefined service standards, measured as a desired number of patients served per facility; \cite{andrade2015abc} define the expected coverage considering both the number of calls received at the call center and the potential increase in the number of ambulances over the planning horizon; \cite{nunez2016column} maximize the sum of some coverage scores that indicate the level of medical access for each truck route, thus optimizing the total level of continuous access for the mobile population; \cite{mousazadeh_accessible_2018} propose a maximal covering objective as a weighted sum of patient flows to the three different medical facilities they consider by multiplying them by their respective coverage scores; and \cite{taymaz2020healthcare} present four different disease-related coverage scores and minimize the so-called weighted expected coverage. Other formulations include participation-based measures (\citealp{dogan_model_2020}), and partial coverage schemes with multiple critical distances and provider collaboration effects (e.g., \citealp{mendoza2024maximal}).

Authors occasionally associate coverage with equity (e.g., \citealp{doerner2007multicriteria}; \citeauthor{smith2009planning}, \citeyear{smith2009planning}, \citeyear{smith2013bicriteria}), or accessibility (e.g., \citealp{murawski2009improving}; \citealp{gunecs2014}). 
In general, objectives involving costs (C), distances (D), or times (T) are connected to accessibility or equity (when they are patient-oriented) or efficiency (when they are system-oriented) criteria. For instance, equity of provision is addressed by minimizing patients' weighted (\citeauthor{smith2009planning}, \citeyear{smith2009planning}, \citeyear{smith2013bicriteria};  \citealp{almeida2021two-step}),   
maximum (\citealp{rouhani2022robust}),  maximum weighted (\citealp{kunkel2014optimal}), and standard deviation of (\citealp{beheshtifar_multiobjective_2015}) travel distances to healthcare facilities, as well as by incorporating constraints to limit maximum travel distances (\citealp{de2022optimization}). \cite{decampos2024multi} define equity in access as minimizing the standard deviation of the coverage ratio across the considered cities.

To enhance accessibility for individuals seeking medical assistance, some works minimize travel costs to nearest facility (\citealp{cocking2012improving}; \citealp{ghaderi2013modeling}; \citealp{kunkel2014optimal}; \citealp{beheshtifar_multiobjective_2015}; \citeauthor{zarrinpoor2018design}, \citeyear{zarrinpoor2017design}, \citeyear{zarrinpoor2018design}; \citealp{chouksey2022hierarchical, chouksey2022optimisation,chouksey2024accelerated};
 \citealp{karakaya2022biobjective}), travel costs upon referrals (\citeauthor{zarrinpoor2018design}, \citeyear{zarrinpoor2017design}, \citeyear{zarrinpoor2018design}; \citealp{chouksey2022hierarchical, chouksey2022optimisation,chouksey2024accelerated}; 
\citealp{karakaya2022biobjective}; \citealp{rouhani2022robust}), treatment costs and accommodation costs for patients' relatives (\citealp{mohammadi2014design}), travel time (\citealp{mohammadi2014design}; \citealp{khodaparasti2017enhancing}; \citealp{mousazadeh2018health}; \citealp{motallebi2020short-term}; \citealp{chen2022multi-objective}; \citealp{rouhani2022robust}), total travel distance (\citealp{galvao2006load}; \citealp{gunecs2014}; \citealp{mendoza-gomez2022location, mendoza-gomez2022regionalization}), or total weighted travel distances (\citealp{yao2014locational}; \citealp{decampos2024multi}; \citealp{mendoza2024regionalization}). This latter objective is treated as an efficiency criterion in the p-median-based models proposed by \cite{smith2013bicriteria}, whereas in their maximal covering-based models, efficiency is defined as the total weighted population served.

Healthcare systems often face the dual challenge of satisfying demand while avoiding congestion. Optimizing the served demand (d) or patients' volume that various medical services can handle is often associated with adjusting service capacity accordingly and consequently improving access to healthcare. Some studies aim to encourage facility placement at high-demand locations by maximizing served demand (\citealp{nunez2016column, de2020roadside, de2022optimization}), whereas others focus on reducing inappropriate or excessive demand to alleviate system costs or congestion (e.g., \citealp{almeida2021two-step, kumar2023location}). Other approaches define accessibility measures that depend on demand, travel time, and capacity of facilities (\citealp{pourrezaie-khaligh2022fix-and-optimize}) or that consider context-specific characteristics of the service environment (\citealp{heyns2024optimisation}). A notable departure from static demand formulations is offered by \cite{alban2022resource}, who model effective adopter demand using a diffusion process.

As mentioned above, system-oriented objectives are often linked with (but not limited to) improving the efficiency of healthcare delivery. Many of these objectives are connected with various types of costs, which might vary in nature. To provide examples, the selected papers consider the minimization of fixed costs associated with the opening, construction, or establishment of new facilities (\citealp{kunkel2014optimal}; \citealp{mohammadi2014design}; \citealp{zahiri2014multi-objective,zahiri2014robust}; \citealp{beheshtifar_multiobjective_2015}; \citealp{shishebori2015robust}; \citealp{mousazadeh2018health}; \citealp{wang2018healthcare}; \citeauthor{zarrinpoor2018design}, \citeyear{zarrinpoor2017design}, \citeyear{zarrinpoor2018design}; \citealp{motallebi2020short-term}; \citealp{chen2022multi-objective}; \citealp{chouksey2022hierarchical, chouksey2024accelerated}; \citealp{mendoza-gomez2022location}; \citealp{rouhani2022robust}) or roads \citep{chen2022multi-objective}, expenses for expansion or improvement of facilities (\citealp{shishebori2015robust}; \citealp{wang2018healthcare}; \citealp{mendoza-gomez2022location};  \citealp{pourrezaie-khaligh2022fix-and-optimize}; \citealp{chouksey2024accelerated}), 
operational costs of facilities (\citealp{ghaderi2013modeling}; \citealp{wang2018healthcare}; \citealp{motallebi2020short-term}; \citealp{pourrezaie-khaligh2022fix-and-optimize}; \citealp{chouksey2024accelerated}), land acquisition costs \citep{beheshtifar_multiobjective_2015}, staffing costs (\citealp{mohammadi2014design}; \citealp{cherkesly2019community}), patient treatment costs (\citealp{mohammadi2014design}; \citealp{zarrinpoor2018design}), transportation costs between different healthcare facilities (\citealp{zahiri2014multi-objective};  \citealp{shishebori2015robust}; \citealp{zarrinpoor2018design}; \citealp{karakaya2022biobjective}; \citealp{rouhani2022robust}), lost service  costs \citep{mohammadi2014design}, costs for exceeding capacity \citep{chouksey2024accelerated}, costs of organ removal or transfer processes (\citealp{zahiri2014multi-objective,zahiri2014robust}; \citealp{rouhani2022robust}), costs of training health workers \citep{vonachen2016optimizing}, equipment and supplies acquisition costs (\citealp{zahiri2014multi-objective, zahiri2014robust}; \citealp{cherkesly2019community}; \citealp{de2022optimization}; \citealp{rouhani2022robust}; \citealp{defreitas2024integrated}), vehicle maintenance costs \citep{cherkesly2019community}, or routing costs (\citealp{vonachen2016optimizing}; \citealp{cherkesly2019community}; \citealp{yang2021outreach}), among others. 
To a lesser extent, system-oriented objectives also involve time- and distance-based measures from the provider's perspective, such as minimizing travel, waiting, or surgery times in transplantation processes (\citealp{rouhani2022robust,zahiri2014multi-objective}), minimizing the time spent on ambulances' relocation or travel (\citealp{andrade2015abc, boutilier2020ambulance}), or minimizing distances from new clinics to non-medical facilities to improve site-suitability (S-s) and compatibility of land-use (\citealp{beheshtifar_multiobjective_2015}). 
Distances can also be associated with improving equity by considering both the maximum distance between different-level facilities offering similar services and the minimum distance between pairs of same-level facilities (\citealp{khodaparasti2017enhancing}). 

Other system-oriented objectives aimed at improving equity in the utilization of health services include minimizing imbalance in facility loading (Im) (\citealp{galvao2006load}), optimizing an envy (En) criterion among demand nodes (\citealp{pourrezaie-khaligh2022fix-and-optimize}), and maximizing the minimum service level (SL) of each patient zone (\citealp{mousazadeh_accessible_2018}).
Furthermore, \cite{mousazadeh_accessible_2018} also address the minimization of system instability (In), defined as the total alterations in assignments and referrals over successive time periods, and \cite{gunecs2014} propose an access measure defined as the expected participation, which decreases with distance. They aim to maximize this measure while also considering equity from both the patients' and physicians' perspectives. For patients, they use the maximum travel distance, and for physicians, they account for the standard deviation of panel sizes (i.e., the number of patients assigned to each physician). Additional objectives include maximizing the number of health workers (nW) (\citealp{smith2009planning,vonachen2016optimizing}), minimizing the number of facilities (nF) to locate (\citealp{eksioglu2024designing}), and maximizing the system quality (Q) (\citealp{wang2018healthcare}) and effectiveness (Effect) (\citealp{de2020roadside}). In these works, quality is associated with ensuring sufficient facility capacity, while effectiveness measures the actual improvement in health outcomes resulting from access to services.
Furthermore, \cite{wang2018healthcare} incorporate patient's environmental benefits (EB) into their model by minimizing the impact of distressing situations (such as congested traffic, noisy surroundings, etc.), thereby fostering a serene medical atmosphere. 
In \cite{chouksey2022hierarchical, chouksey2022optimisation}, availability (Av) aspects are taken into consideration, which depend on whether there are referral facilities within the maximum coverage distance of the lower-level facilities and whether the total capacity of the healthcare system meets the demand. Instead of being incorporated as objectives, these high-level criteria are treated as constraints, typically by imposing minimum or maximum threshold requirements. \cite{chouksey2024accelerated}, on the other hand, address availability by incorporating facility capacity limits and a penalty for exceeding them directly into the objective function.

Overall, the literature shows a clear trend toward considering multiple types of objectives to address accessibility, equity, and efficiency. While the most common approaches focus on optimizing coverage and distance-based measures from the patient side and minimizing costs from the system side, several other objectives have also been considered, as discussed above. Many of the selected papers examine various models with distinct objectives, which they solve individually. While certain studies adopt a unified objective that encompasses various similar features (such as the sum of different costs), others address conflicting objectives using a multi-objective function, as we will elaborate on in Section~\ref{subsec:obj_function}.

\subsubsection{Hierarchy}\label{subsec:hierarchy}
In Section~\ref{subsec:type_of_care}, we observed that the selected papers predominantly focus on a single level of care, whether primary, secondary, or tertiary. While it is common for a particular level of care or treatment to be delivered in a specific type of facility, it is evident that patients often necessitate various services, frequently requiring specialized care subsequent to initial treatment. Healthcare systems typically present a hierarchical structure, wherein the level of care is determined by the nature of the required treatment.

For instance, Malaysia's public healthcare system offers three levels of service (primary, secondary, and tertiary), provided across four types of facilities (community clinics, district health offices and hospitals, state health departments, and specialized hospitals and medical college hospitals) operating at different administrative levels (community village, district, state, and national). However, \cite{shariff2012location} exclusively study community clinics, which offer primary-level services, including health clinics and rural clinics. Similarly, \cite{kumar2023location}, while acknowledging the hierarchical nature of the Indian system, does not incorporate it into their proposed location models. Nevertheless, the author considers two types of primary care centers, aiming to establish them in a manner that diverts patients requiring primary services away from specialized hospitals, thus relieving the overcrowding of their outpatient clinics. This issue is also related to the referral systems in developing settings, where, unlike in developed countries, patients may sometimes access higher-level facilities without needing referrals from lower-level ones. 
Vaccination supply chains in developing countries also present hierarchical structures, which generally include a storage central location at the state level, distribution hubs at regional or district levels, and session sites where vaccines are actually administered. However, the studies addressing vaccination supply chains (\citealp{lim2016coverage}; \citealp{srivastava2021strengthening}; \citealp{yang2021outreach}; \citealp{eksioglu2024designing}) do not consider hierarchical location models but focus their efforts on strengthening the last tier of the network.  

Some of the selected articles (14 out of 56) do integrate this hierarchy when formulating their problems. For instance, \cite{smith2009planning} propose a series of hierarchical models tailored for the planning of community health programs, where service provision may occur through village health workers, community facilities, or hospitals. The exploration of different referral strategies between these tiers of healthcare provision is also considered. 
\cite{smith2013bicriteria} also introduce several hierarchical models for a general health system and account for different types of hierarchies, referrals, flow disciplines, and access modes to the system.
Along the lines of the former work, \cite{khodaparasti2017enhancing} introduce a multi-service hierarchical location model to enhance community health schemes in Iran. There are three types of service providers (community-based organizations, consulting centers, and addiction treatment clinics) that adhere to a distinctly tiered hierarchy, where upper-level facilities exclusively provide services corresponding to their own tier. Both consulting centers and addiction treatment clinics assist their own patients, with only a portion of their users being referred by community-based organizations. Moreover, the authors incorporate the preferences of referred recipients into their model. These preferences, which are related to their proximity to second-level facilities, are integrated to reinforce patients' commitment and adherence to treatment.

Although \cite{zahiri2014multi-objective,zahiri2014robust} and \cite{rouhani2022robust} analyze the Iranian organ transplant supply chain, only the latter study provides a location model that considers the inherent hierarchy of this network, where the three levels of facilities involved are  hospitals, organ procurement units, and transplant centers. In addition to considering different facility levels, these authors also propose a two-tier hierarchical method for organ allocation that includes regional zones and the entire country. Meanwhile, \cite{chouksey2022hierarchical,chouksey2022optimisation,chouksey2024accelerated} incorporate the hierarchical structure of the Indian maternal healthcare system into their models, which consists of three types of interacting facilities that are successively inclusive. In particular, low-level facilities are primary centers that offer basic services, mid-level facilities offer also Cesarean procedures, and high-level facilities can also assist neonatal or unplanned deliveries. Mothers-to-be can be referred from any lower-level facility to an upper-level one, even to a higher-level service offered at the same facility. \cite{galvao2006load} study a similar hierarchical maternal healthcare system in the Brazilian context, where referrals are restricted to occur only from mid-level to highest-level facilities. 
In contrast, \cite{karakaya2022biobjective} explicitly allow backward referrals of mothers and their newborns from highest- to mid-level facilities in the Turkish maternal healthcare system.

\cite{mousazadeh_accessible_2018, mousazadeh2018health}, \citeauthor{zarrinpoor2018design} (\citeyear{zarrinpoor2017design}, \citeyear{zarrinpoor2018design}), and \cite{motallebi2020short-term} propose hierarchical location-allocation models to capture the intricacies of the Iranian healthcare system. In particular, \cite{mousazadeh_accessible_2018, mousazadeh2018health} consider three types of facilities: primary healthcare centers offering basic and preventive health services, clinics or regional healthcare centers providing more invasive care, and hospitals delivering specialized services; \citeauthor{zarrinpoor2018design} (\citeyear{zarrinpoor2017design},\citeyear{zarrinpoor2018design}) examine a two-level multi-flow hierarchy, where lower-level hospitals (district hospitals) deliver community-based and non-specialized services, while higher-level hospitals (regional hospitals or central hospitals) offer both community-based and specialized services; and \cite{motallebi2020short-term} introduce a bi-hierarchy model that distinguishes between clinics offering outpatient services and hospitals providing also inpatient and emergency services. 
Furthermore, although all these models account for referrals, it is noteworthy that such referrals are not mandatory for patients to receive care at higher-level facilities. In Iran, as well as in other developing countries, individuals have the option to directly access services at higher-level centers.

Other hierarchical-based strategies have been considered in the literature. For instance, \cite{wang2018healthcare} develop a bi-level model in which the upper-level layer incorporates location decisions for healthcare centers. Despite describing the distinct hierarchies of China's healthcare system in rural regions, the authors do not account for varying facility levels. In the lower-level layer, capacity decisions are determined based on optimal locations.  
In a different context, \cite{heyns2024optimisation} consider improving healthcare access by building road sequences that result from grouping individual roads, which are classified into three connectivity hierarchies based on the existing road network. Tier 1 roads are feasible for construction independently, whereas tiers 2 and 3 require the prior construction of lower-level tiers.

\subsubsection{Uncertainty}\label{subsec:uncertainty}
Strategic decisions are often hindered by a multitude of uncertainties stemming from various sources, such as demographic factors, resource constraints, and unpredictable demand patterns. In developing countries, these uncertainties are even greater due to factors like lack of data records, weakness of healthcare systems often linked to political instability, or infrastructural vulnerability. Addressing these uncertainties becomes crucial for policymakers and healthcare administrators to design robust and effective strategies that can withstand the fluctuations inherent in such environments.

Despite the great importance of this matter, only 19 of the selected papers consider non-deterministic parameters. These include demands (e.g., \citealp{mohammadi2014design}; \citealp{zahiri2014multi-objective}; \citealp{shishebori2015robust}; \citealp{boutilier2020ambulance}; \citealp{yang2021outreach}), service times (e.g., \citealp{zahiri2014multi-objective}; \citealp{andrade2015abc};  \citealp{dogan_model_2020};  \citealp{motallebi2020short-term}), system costs (e.g., \citealp{zahiri2014multi-objective}; \citealp{mousazadeh_accessible_2018}; \citealp{wang2018healthcare}), travel times (e.g., \citealp{mohammadi2014design};  \citealp{zahiri2014multi-objective};  \citealp{boutilier2020ambulance}; \citealp{yang2021outreach}), service capacity (\citealp{mohammadi2014design}; \citeauthor{zarrinpoor2018design}, \citeyear{zarrinpoor2017design}, \citeyear{zarrinpoor2018design}), availability of resources (\citealp{zahiri2014multi-objective}; \citealp{mousazadeh_accessible_2018}), or supplies fluctuations (\citealp{zahiri2014multi-objective}; \citealp{rouhani2022robust}), among others. Problem-specific parameters such as the time for organ removal (\citealp{zahiri2014multi-objective}) and the number of donors (\citealp{zahiri2014robust}) for organ transplant networks, as well as disease incidence (\citealp{wang2018healthcare}), have also been considered as uncertain to a lesser extent. In addition, \citeauthor{zarrinpoor2018design} (\citeyear{zarrinpoor2017design}, \citeyear{zarrinpoor2018design}) account for stochasticity in system reliability, which is the system's ability to maintain service when one or more facilities are disrupted. The latter study also considers uncertainty in referral rates and geographical accessibility, this latter measured as the maximum acceptable distance to access healthcare facilities. The covering distance is also treated as uncertain in \cite{mohammadi2014design}. \cite{mousazadeh_accessible_2018} is the only study to adopt soft constraints and uncertainty in the target values of goals, which is coped with by a flexible programming approach.

Unsurprisingly, demands emerge as the predominant element when it comes to contemplating stochastic behavior, with all but one paper addressing this aspect. Researchers have harnessed various methodologies to account for demands' uncertainties, such as fuzzy numbers (\citealp{zahiri2014multi-objective}), queuing theory (\citealp{mohammadi2014design}; \citealp{andrade2015abc}; \citeauthor{zarrinpoor2018design}, \citeyear{zarrinpoor2017design}, \citeyear{zarrinpoor2018design};  \citealp{dogan_model_2020}; \citealp{motallebi2020short-term}), stochastic programming (\citealp{wang2018healthcare}; \citealp{yang2021outreach}), robust optimization (\citealp{shishebori2015robust}; \citealp{zarrinpoor2017design}; \citealp{boutilier2020ambulance}; \citealp{motallebi2020short-term}; \citealp{rouhani2022robust}), or scenario-based approaches (\citealp{shishebori2015robust}; \citealp{zarrinpoor2018design}; \citealp{boutilier2020ambulance}; \citealp{taymaz2020healthcare}), among others. \cite{alban2022resource} model the effective adopter demand using the Bass diffusion model  \citep{bass1969}, and estimate its parameters through a machine-learning approach based on gradient boosting. To address uncertainty in these parameter estimates, the authors propose two dynamic heuristics that adapt the allocation of visit capacities over time. In doing so, they also account indirectly for demand uncertainty, since demand is a function of the estimated parameters. Uncertainty on service times is always accounted for by employing queuing theory. Other methodologies in the selected literature include interval programming or interval-valued fuzzy set (\citealp{mohammadi2014design}), robust possibilistic programming (\citealp{zahiri2014robust}; \citealp{mousazadeh2018health}), chance-constrained approaches (\citealp{khodaparasti2017enhancing}), or hybrid optimization and simulation frameworks (\citealp{chouksey2022optimisation}). 

Table~\ref{tab:uncertainty} expands upon the information discussed herein, outlining the uncertain elements and the corresponding methods chosen to address them for each of the 19 papers.

\begin{table}[ht!]\footnotesize
\caption{Uncertain parameters and methods to address the uncertainty.}\label{tab:uncertainty}
\resizebox{\linewidth}{!}{
\begin{tabular}{p{0.5cm}p{5.8cm}ll}
\noalign{\smallskip}\hline\noalign{\smallskip} 
No. & Reference & Parameter & Method  \\
\noalign{\smallskip}\hline\noalign{\smallskip}

    11 & \cite{mohammadi2014design}     & Demands, service time & Queuing theory \\
    & & Travel time & Scenario-based approach 
    \\
    & & Service capacity, system costs & Interval programming\\
    & & Covering distance & Interval-valued fuzzy set \\[0.2em] \arrayrulecolor{gray!30}\hline

    13 & \cite{zahiri2014multi-objective}      & \multirow{3}{*}{\parbox{6cm}{Resources, service time
    \\ System costs, travel time, supplies, \\demands, time for removal process}} 
    & \multirow{3}{*}{\parbox{3cm}{Queuing theory \\ Fuzzy numbers  
    \\  }}  \\
    & & & \\
     & & & \\[0.2em] \arrayrulecolor{gray!30}\hline

    14 & \cite{zahiri2014robust}   &  \multirow{2}{*}{\parbox{6cm}{System costs, transportation\\ times, number of donors, demands}} &   \multirow{2}{*}{\parbox{4cm}{Robust possibilistic programming\\}}\\
    & & & \\[0.2em] \arrayrulecolor{gray!30}\hline

    15 & \cite{andrade2015abc}     &   Demands, service time  &  Queuing theory \\[0.2em] \arrayrulecolor{gray!30}\hline

    17 & \cite{shishebori2015robust}      &  Demands, transfer costs & \multirow{2}{*}{\parbox{6cm}{Scenario-based approach\\ and robust optimization}} \\
    & & & \\[0.2em] \arrayrulecolor{gray!30}\hline

    22 & \cite{zarrinpoor2017design}    & Demands, service capacity &  Queuing theory
    \\
    & & System reliability, demand rates & Two-stage robust optimization
    \\[0.2em] \arrayrulecolor{gray!30}\hline

    23 & \cite{khodaparasti2018multi-period}    & \multirow{2}{*}{\parbox{4cm}{Demands\\}} & \multirow{2}{*}{\parbox{5cm}{Distributionally robust\\ chance-constrained approach}} 
    \\
    & & & \\[0.2em] \arrayrulecolor{gray!30}\hline

    24 & \cite{mousazadeh_accessible_2018}    & Resources, system costs &  Fuzzy numbers
    \\
    & & Soft constraints, target values of goals & Flexible programming approach
    \\[0.2em] \arrayrulecolor{gray!30}\hline

    25 & \cite{mousazadeh2018health}    & Demands, opening costs,

    &  Robust possibilistic programming
    \\
    & & service capacity, transportation times 
    & 
    \\[0.2em] \arrayrulecolor{gray!30}\hline

    26 & \cite{wang2018healthcare}    & Demands, disease incidence  & Stochastic programming 
    \\
    & & System costs 
    & Fuzzy numbers\\[0.2em] \arrayrulecolor{gray!30}\hline
    
    27 & \cite{zarrinpoor2018design}    &   \multirow{3}{*}{\parbox{6cm}{Service capacity, system reliability,\\ demands and referral rates,\\ geographical accessibility}}  &  \multirow{3}{*}{\parbox{5cm}{Queuing theory and robust\\ scenario-based approach}} 
    \\
    & & & \\
    & & & \\[0.2em] \arrayrulecolor{gray!30}\hline

    29 & \cite{boutilier2020ambulance}   &  Travel time, demands & Scenario-based approach and robust optimization
    \\[0.2em] \arrayrulecolor{gray!30}\hline

    31 & \cite{dogan_model_2020}    & Demands, service time  & Queuing theory \\[0.2em] \arrayrulecolor{gray!30}\hline

    32 & \cite{motallebi2020short-term}    & Demands, service time & Queuing theory \\
    & & Demand rates & Robust optimization\\[0.2em] \arrayrulecolor{gray!30}\hline

    33 & \cite{taymaz2020healthcare}     &  Demands & Scenario-based approach
    \\[0.2em] \arrayrulecolor{gray!30}\hline

    37 & \cite{yang2021outreach}   &  Demands, travel time & Stochastic programming \\[0.2em] \arrayrulecolor{gray!30}\hline

    38 & \cite{alban2022resource}    & Bass model parameters &  Dynamic heuristics
    \\[0.2em] \arrayrulecolor{gray!30}\hline

     41 & \cite{chouksey2022optimisation}    &  \multirow{2}{*}{\parbox{4cm}{Demands\\}} & \multirow{2}{*}{\parbox{4cm}{Hybrid optimization and \\simulation framework}}  \\
    & & &\\[0.2em] \arrayrulecolor{gray!30}\hline

    48 & \cite{rouhani2022robust}     &  Demands, supplies &  Robust optimization  \\

     \arrayrulecolor{black}
\noalign{\smallskip}\hline\noalign{\smallskip}
    \end{tabular}
    }
\end{table}

\subsubsection{Multi-period settings}
\label{subsec:multiperiod}
The majority of the studies in the surveyed literature adopt a single-period planning perspective, modeling healthcare accessibility decisions as static and assuming that demand patterns and system characteristics remain constant over time. While appropriate for certain strategic analyses, this approach may overlook important 
temporal dynamics such as demographic changes, gradual infrastructure development, the progressive availability of resources, or phased policy implementation. Multi-period planning models address these limitations by explicitly incorporating time as a modeling dimension, allowing decisions and system parameters to evolve across the planning horizon \citep{nickel2019}. Despite their relevance, only 12 out of the 56 selected papers consider multi-period settings, as showcased in the tenth column of Table~\ref{tab:problem_features}. 

The reviewed papers incorporate multi-period settings in different ways. A first group allows strategic decisions, such as the opening of facilities or network links, to occur at different time periods, typically reflecting budgetary or resource constraints (e.g., \citealp{ghaderi2013modeling, khodaparasti2018multi-period, pourrezaie-khaligh2022fix-and-optimize, chouksey2024accelerated}). Most multi-period studies, however, adopt formulations in which strategic decisions remain fixed over the planning horizon, while operational decisions are updated over time. Examples include models where outreach locations are determined once, but routing and service assignments are periodically revised (e.g., \citealp{yang2021outreach}); applications involving the allocation of emergency vehicles to stand-by-points over a planning horizon (e.g., \citealp{andrade2015abc}); and models with time-varying demand or supplies (e.g., \citealp{zahiri2014multi-objective, zahiri2014robust, dogan_model_2020, rouhani2022robust, eksioglu2024designing}).


\subsubsection{Geographical area and country of application}
\label{subsec:geography}
The last two columns of Table~\ref{tab:problem_features} indicate the geographical scale at which the case study is conducted and the country of application. In addition, they highlight whether the areas considered are predominantly rural (R), predominantly urban (U), or a combination of both (R, U). This classification is based on the descriptions provided by the authors of the original studies regarding the characteristics of the analyzed regions, as well as a close inspection of regions, districts, or provinces using open-source databases such as \textit{City Population} \citep{Brinkhoff}.

We note that the majority of the selected papers concentrate their case studies on one or more administrative units of the same category. For instance, \cite{yao2014locational} and \cite{vonachen2016optimizing} examine various rural districts in Mozambique and Liberia, respectively, while \cite{shishebori2015robust} investigate multiple towns in Iran. \cite{srivastava2021strengthening} explore eight Indian districts and \cite{mendoza-gomez2022location} study
17 municipalities in the northern zone of the State of Mexico, both encompassing rural and urban areas. Only a few works address two different types of administrative units. \cite{shariff2012location} center their case study on a district with populations residing in both rural and urban settings, as well as a rural sub-district in Malaysia. \cite{kunkel2014optimal} initially concentrate on a rural district in Malawi before extending their application to the national level. \cite{zahiri2014multi-objective} and \cite{de2022optimization} 
similarly investigate the entire country for their case studies.  \cite{nunez2016column}, \cite{de2020roadside}, and \cite{taymaz2020healthcare} consider the major transport corridors across multiple countries in South-East or Sub-Saharan Africa, and \cite{yang2021outreach} choose various rural regions within four of these countries, also without specifying them. \cite{alban2022resource} focus their case study in hard-to-reach, rural sites across Uganda. Additional considerations on the country of application are available in the Supplementary Material.


\subsection{Modeling and solution approaches}\label{subsec:opt_approach}
Next, we focus on the modeling and solution approaches for the studied papers, considering the type of problems analyzed, the mathematical models selected for their formulation, the type of objective functions employed, and the proposed solution methods. Table~\ref{tab:opt_approach} provides a summary of these aspects.

\subsubsection{Type of problem}\label{subsec:type_of_problem}
In this section, we categorize the selected papers based on the type of optimization problem they address. These include location problem (LP), allocation problem (AP), location-routing problem (LRP), network design problem (NDP), and location network design problem (LNDP). LP includes studies that focus on location, with or without explicit allocation modeling. NDP includes studies that focus solely on transportation network improvements, while LNDP refers to studies integrating both facility location and network design decisions.

\scriptsize
\setlength\LTleft{-1cm}
\setlength\LTright{-1cm}

\begin{longtable}{  
    p{0.4cm}  
    p{4.9cm} 
    !{\color{gray!30}\vrule} 
   l 
    !{\color{gray!30}\vrule}
    p{3cm} 
    !{\color{gray!30}\vrule}
    p{1.7cm} 
    !{\color{gray!30}\vrule}
    l 
    !{\color{gray!30}\vrule}
    p{2.7cm}  
    }
\caption{Modeling and solution approaches for the selected papers.\label{tab:opt_approach}}\\ 
\noalign{\smallskip}\hline\noalign{\smallskip} 
& & \multirow{2}{*}{\parbox{1.5cm}{Type of problem}} & \multirow{2}{*}{\parbox{3cm}{Mathematical model}} & \multirow{2}{*}{\parbox{2.6cm}{Objective\\ function}} &  \multicolumn{2}{l}{\multirow{2}{*}{\parbox{5cm}{\centering Solution approach}}} \\
No. & Reference & & & & & \\
&  & & & & Exact & Approximate \\
\noalign{\smallskip}\hline\noalign{\smallskip}
\endfirsthead

\caption*{Table~\ref{tab:opt_approach} (continued): Modeling and solution approaches for the selected papers.}\\
\noalign{\smallskip}\hline\noalign{\smallskip} 
& & \multirow{2}{*}{\parbox{1.5cm}{Type of problem}} & \multirow{2}{*}{\parbox{3cm}{Mathematical model}} & \multirow{2}{*}{\parbox{2.6cm}{Objective\\ function}} &  \multicolumn{2}{l}{\multirow{2}{*}{\parbox{5cm}{\centering Solution approach}}} \\
No. & Reference & & & & & \\
&  & & & & Exact & Approximate \\
\noalign{\smallskip}\hline\noalign{\smallskip}
\endhead

\noalign{\smallskip}\hline\noalign{\smallskip}
\endfoot

\noalign{\smallskip}\hline\noalign{\smallskip}
\multicolumn{7}{l}{\parbox{1.05\linewidth}{\textbf{Type of problem}: LP (location problem), LRP (location-routing problem), NDP (network design problem), LNDP (location network design problem), AP (allocation problem); 
\textbf{Mathematical model}: MILP (mixed-integer linear programming), MINLP (mixed-integer nonlinear programming), ILP (integer linear programming), FP (fuzzy programming), SP (stochastic programming), SDP (stochastic dynamic programming), RO (robust optimization), FGP (fuzzy goal programming), INLP (integer nonlinear programming), BLP (bi-level programming), GP (goal programming), NLP (nonlinear programming); 
\textbf{Objective function}: S (single-objective), M/no./method (multi-objective/no. of objectives/method to cope with the multi-objective problem), $\varepsilon$ ($\varepsilon$-constraint method), EA (evolutionary algorithm), WS (weighted sum), GT (game theory), TH \citep{Torabi2008}, FGP (fuzzy goal programming), GP (goal programming), SA (simulated annealing);
\textbf{Solution approach}: CS (commercial solver), H (heuristics), Lag (Lagrangian), MH (metaheuristics), AC (ant colony), GA (genetic algorithm), Gr (greedy), FO (fix-and-optimize), SA (simulated annealing), ICA (imperialist competitive algorithm), HSA (hybrid solution approach), DE (differential evolution), ABC (artificial bee colony), DP (decision process), TOPSIS (technique for order of preference by similarity to ideal solutions), CG (column generation), BD (Benders decomposition), OA (outer approximation), PSO (particle swarm optimization), VG (variable generation), RCG (row and column generation), SG (scenario generation), O (other), RF (relax-and-fix), TDH (top-down heuristic), EA (evolutionary algorithm), IGA (iterated greedy algorithm).}} \\
\endlastfoot

    1 & \cite{galvao2006load} & LP & MILP, MINLP (corr. MILP) & S, M/2/$\varepsilon$ & CS (CPLEX) & H (Lag) \\[0.2em] \arrayrulecolor{gray!30}\hline

    2 & \cite{doerner2007multicriteria} & LRP & ILP & M/3/EA & & MH (AC, GA) \\[0.2em] \arrayrulecolor{gray!30}\hline

    3 & \cite{murawski2009improving}     &   NDP     &  ILP &  S  &  CS (CPLEX) &  \\[0.2em]
     \arrayrulecolor{gray!30}\hline

    4 & \cite{smith2009planning}     & LP  &  MILP    &  S &  CS (Xpress-MP)
    & \\[0.2em] \arrayrulecolor{gray!30}\hline
    
    5 & \cite{cocking2012improving}       & LNDP    &  MILP & S  & CS (CPLEX)  & \\[0.2em] \arrayrulecolor{gray!30}\hline
    
    6 & \cite{shariff2012location}       &  LP  &  MILP  &  S   &  CS (CPLEX) &  MH (GA) \\[0.2em] \arrayrulecolor{gray!30}\hline
    
    7 & \cite{ghaderi2013modeling}     &   LNDP  &  MINLP (corr. MILP) & S   &  CS (CPLEX) & H (Gr), MH (FO) 
    \\[0.2em] \arrayrulecolor{gray!30}\hline

    8 & \cite{smith2013bicriteria} & LP & MILP & M/2/WS & CS (Xpress-MP) &  \\[0.2em] \arrayrulecolor{gray!30}\hline

    9  & \cite{gunecs2014}    & LP   & MILP &  S  & CS 
    &  
    \\[0.2em] \arrayrulecolor{gray!30}\hline

    10 & \cite{kunkel2014optimal}   &  LP   & ILP  & S   &  CS (Gurobi)   & \\[0.2em] \arrayrulecolor{gray!30}\hline

    11 & \cite{mohammadi2014design}      &  LNDP   &  FP, SP 
    &   M/2/GT
    &  CS (BARON) & MH (SA, ICA), HSA \\[0.2em] \arrayrulecolor{gray!30}\hline
    
    12 & \cite{yao2014locational}      &   LP  &  ILP & S   & CS (Gurobi)  & \\[0.2em] \arrayrulecolor{gray!30}\hline
    
    13 & \cite{zahiri2014multi-objective}      &  LP   &  \parbox{3.2cm}{MINLP (corr. MILP),\\ FP
     } &   M/2/TH
    & CS (BARON)  & MH (SA, DE)\\[0.2em] \arrayrulecolor{gray!30}\hline

    14 & \cite{zahiri2014robust}      &  LP   & MILP, FP   & S   &  CS  & \\[0.2em] \arrayrulecolor{gray!30}\hline

    15 & \cite{andrade2015abc}     &  LP   & SDP
    &  S    & CS (CPLEX)  & MH (ABC)\\[0.2em] \arrayrulecolor{gray!30}\hline

    16 & \cite{beheshtifar_multiobjective_2015}    & LP  &  MINLP &    M/4/EA &  & MH (GA), DP (TOPSIS)\\[0.2em] \arrayrulecolor{gray!30}\hline

    17 & \cite{shishebori2015robust}      & LNDP    &  MILP, RO  &  S  & CS (CPLEX) &  \\[0.2em] \arrayrulecolor{gray!30}\hline

    18 & \cite{lim2016coverage}    & LP  &  MILP  &  S   &  CS (CPLEX) &\\[0.2em] \arrayrulecolor{gray!30}\hline

    19 & \cite{nunez2016column}    &  LP   & MILP  & M/3/WS    &   & H (CG)\\[0.2em] \arrayrulecolor{gray!30}\hline
    
    20 & \cite{vonachen2016optimizing}    &  LRP   &  ILP &  S   & CS (CPLEX) & \\[0.2em] \arrayrulecolor{gray!30}\hline

    21 & \cite{khodaparasti2017enhancing}  &  LP   &  MILP, FGP    & M/3/FGP & CS (CPLEX)  & \\[0.2em] \arrayrulecolor{gray!30}\hline

    22 & \cite{zarrinpoor2017design} & LP & MILP, RO & S & CS (CPLEX), BD &  \\[0.2em] \arrayrulecolor{gray!30}\hline
    
    23 & \cite{khodaparasti2018multi-period}     &  LP   &  \parbox{3.2cm}{SP (eq. INLP w/corr. ILP)}
    &  S & OA & \\[0.2em] \arrayrulecolor{gray!30}\hline

    24 & \cite{mousazadeh_accessible_2018}    &  LP   & \parbox{3.2cm}{MINLP (corr. MILP),\\ FP}  
    &    M/4/$\varepsilon$ & CS (CPLEX) &  \\[0.2em] \arrayrulecolor{gray!30}\hline

    25 & \cite{mousazadeh2018health} & LP & MINLP (corr. MILP), FP & M/2/TH & CS (CPLEX) &  \\[0.2em] \arrayrulecolor{gray!30}\hline
    
    26 & \cite{wang2018healthcare}    & LP &  
    BLP, SP
    &   M/4/EA
    &   & MH (PSO)\\[0.2em] \arrayrulecolor{gray!30}\hline
    
    27 & \cite{zarrinpoor2018design}    &   LP  &  \parbox{3.2cm}{ MINLP (corr. convex model),
    \\ RO}

    & S   & CS (CPLEX), BD & \\[0.2em] \arrayrulecolor{gray!30}\hline

    28 & \cite{cherkesly2019community}   &  LRP   & MILP  & S   & CS (CPLEX), VG & H \\[0.2em] \arrayrulecolor{gray!30}\hline
    
    29 & \cite{boutilier2020ambulance}   &  LRP   & RO (eq. MILP) 
    &S & RCG  &  H (SG)  \\[0.2em] \arrayrulecolor{gray!30}\hline

    30 & \cite{de2020roadside}    & LP  &  MINLP (corr. MILP)  &  M/2/WS &  CS (CPLEX) &    \\[0.2em] \arrayrulecolor{gray!30}\hline

    31 & \cite{dogan_model_2020}    & LP  &   MILP, GP  &   M/3/GP   & CS (CPLEX) &\\[0.2em] \arrayrulecolor{gray!30}\hline
        
    32 & \cite{motallebi2020short-term}    &   LP  & MINLP (corr. MILP), RO  &   M/2/EA

     & CS (CPLEX) & MH (GA)\\[0.2em] \arrayrulecolor{gray!30}\hline

    33 & \cite{taymaz2020healthcare}     &  LP   &  SP  &  S   & CS (CPLEX) & \\[0.2em] \arrayrulecolor{gray!30}\hline
    
    34 & \cite{almeida2021two-step}    &   LP   &  MILP &  S  & CS (GLPK) & \\[0.2em] \arrayrulecolor{gray!30}\hline

    35 & \cite{flores2021optimizing}    &   LP    &  ILP  &  S  & O & \\[0.2em] \arrayrulecolor{gray!30}\hline
    
    36 & \cite{srivastava2021strengthening}    &  LP   &  ILP  & S  & CS (Excel Solver) & H (Gr) \\[0.2em] \arrayrulecolor{gray!30}\hline

    37 & \cite{yang2021outreach}    &  LRP & MILP, SP
    &  S     & CS (Gurobi)  & \\[0.2em] \arrayrulecolor{gray!30}\hline
    
   38 & \cite{alban2022resource} & AP & NLP & S & & H (Lag, O) \\[0.2em] \arrayrulecolor{gray!30}\hline

    39 & \cite{chen2022multi-objective}   &  LNDP   & ILP  &   M/3/SA
    &  
    & MH (SA)\\[0.2em] \arrayrulecolor{gray!30}\hline

    40 & \cite{chouksey2022hierarchical}   &   LP  & MILP  &  S   &  CS (Gurobi) &  H (RF+FO)  \\[0.2em] \arrayrulecolor{gray!30}\hline

    41 & \cite{chouksey2022optimisation}    &   LP  &  MILP & S   & CS (Gurobi) & \\[0.2em] \arrayrulecolor{gray!30}\hline

   42 & \cite{de2022optimization}    &  LP    & ILP  &   S  & CS (GLPK) & \\[0.2em] \arrayrulecolor{gray!30}\hline

    43 & \cite{elorza2022assessing}    &   LP
    &  MILP &  S     &  CS (CPLEX) &\\[0.2em] \arrayrulecolor{gray!30}\hline

    44 & \cite{karakaya2022biobjective} & LP & MILP & M/2/WS & CS (CPLEX) & H (TDH, Lag) \\[0.2em] \arrayrulecolor{gray!30}\hline

    45 & \cite{mendoza-gomez2022location}      &  LP   & ILP   &  S, M/2/$\varepsilon$ & CS (CPLEX) &  \\[0.2em] \arrayrulecolor{gray!30}\hline

    46 & \cite{mendoza-gomez2022regionalization}   & LP  & ILP  & S  &  CS (CPLEX) &    \\[0.2em] \arrayrulecolor{gray!30}\hline
    
    47 & \cite{pourrezaie-khaligh2022fix-and-optimize} &  LNDP  &  MINLP (corr. MILP) & 
      M/3/$\varepsilon$  &    &  MH (FO)
        \\[0.2em] \arrayrulecolor{gray!30}\hline

    48 & \cite{rouhani2022robust}       &  LP   &  \parbox{3.2cm}{MINLP (corr. MILP),\\ RO} &   M/3/$\varepsilon$ &  CS  & \\[0.2em] \arrayrulecolor{gray!30}\hline

    49 & \cite{kumar2023location}     &  LP  &  ILP, MILP &  S
    &  CS (CPLEX) & H \\[0.2em] \arrayrulecolor{gray!30}\hline

    50 & \cite{decampos2024multi} & LP & MINLP & M/3/EA & & MH (GA, EA) \\[0.2em] \arrayrulecolor{gray!30}\hline

    51 & \cite{defreitas2024integrated} & LP & MINLP (corr. MILP) & S & CS (GLPK) &  \\[0.2em] \arrayrulecolor{gray!30}\hline

    52 & \cite{eksioglu2024designing}   &   LP  &  ILP, MILP &  M/2/WS & CS (CPLEX), BD &   \\[0.2em] \arrayrulecolor{gray!30}\hline

    53 & \cite{heyns2024optimisation}    &  NDP  &  ILP &  M/3/WS, M/3/EA &  CS (CPLEX) & MH (GA)  \\[0.2em] \arrayrulecolor{gray!30}\hline

    54 & \cite{mendoza2024maximal} & LP & ILP & S & CS (CPLEX) &  \\[0.2em] \arrayrulecolor{gray!30}\hline

    55 & \cite{mendoza2024regionalization} & LP & ILP & S & CS (CPLEX) & MH (IGA) \\[0.2em] \arrayrulecolor{gray!30}\hline

    56 & \cite{chouksey2024accelerated} & LP & MILP & S & CS (Gurobi), BD & H, MH(FO + SA)
     \\
\end{longtable}



\normalsize

A substantial share of the reviewed papers (42 out of 56) formulates healthcare access as location or location-allocation problems, typically building on classical discrete location models. In particular, p-median formulations \citep{hakimi1964} are frequently used to minimize travel distances when locating health facilities, services, or workers (e.g., \citealp{smith2009planning, yao2014locational}), while maximal covering models \citep{church1974} are widely adopted when coverage within a critical distance is the primary objective (e.g., \citealp{smith2009planning, lim2016coverage, elorza2022assessing, mendoza2024maximal}). Set covering formulations \citep{toregas1973} appear when full coverage is required (e.g., \citealp{eksioglu2024designing}), and more recent extensions, such as the cooperative maximal covering model of \cite{li2018} have been used to capture interactions among facilities (e.g., \citealp{flores2021optimizing}).

Beyond pure location decisions, most studies explicitly integrate allocation decisions: 34 out of the 56 selected papers jointly determine where facilities or services should be located and how patients or healthcare resources should be assigned.
To name a few, \cite{shariff2012location} and \cite{beheshtifar_multiobjective_2015} attempt to determine the optimum placement of new facilities and their demand allocation using the capacitated maximal covering location problem \citep{current1988} and the p-median approach, respectively. Other works seeking to simultaneously locate health centers and allocate patients are \cite{galvao2006load}, \cite{smith2013bicriteria}, \cite{gunecs2014}, \cite{khodaparasti2017enhancing, khodaparasti2018multi-period}, \cite{mousazadeh_accessible_2018,mousazadeh2018health}, \citeauthor{zarrinpoor2018design} (\citeyear{zarrinpoor2017design},\citeyear{zarrinpoor2018design}), \cite{chouksey2022hierarchical, chouksey2022optimisation,chouksey2024accelerated}, \cite{karakaya2022biobjective}, \cite{mendoza-gomez2022location, mendoza-gomez2022regionalization, mendoza2024regionalization}, \cite{rouhani2022robust}, and \cite{kumar2023location}. \cite{mendoza-gomez2022location, mendoza-gomez2022regionalization,mendoza2024regionalization} and \cite{chouksey2024accelerated} also attempt to upgrade existing facilities by expanding their capacity or providing additional healthcare services, while \cite{gunecs2014} explicitly model the allocation of physicians to facilities. \cite{zahiri2014multi-objective,zahiri2014robust}, besides the optimal location of hospitals and transplant centers, also consider the allocation of appropriate equipment or shipping agents. 
\cite{kunkel2014optimal} employ the p-median model to assign health workers to geographical units and the capacitated facility location problem to assign backpacks to workers and the latter to resupply centers. \cite{andrade2015abc} analyze how to optimally determine ambulance stations' locations and the allocation and repositioning of these emergency care vehicles to so-called ``stand-by points''. \cite{wang2018healthcare} propose a bi-level approach where the upper level focuses on location decisions and the lower level addresses the capacity allocation. \cite{mousazadeh2018health},
\cite{de2020roadside}, \cite{motallebi2020short-term},  \cite{taymaz2020healthcare}, and \cite{almeida2021two-step} aim to locate new healthcare facilities and to decide which services each should provide. In particular, \cite{mousazadeh2018health} determine the optimal specialty service portfolio at the third level of the hierarchy; \cite{motallebi2020short-term} model the allocation of service units to the facilities; \cite{de2020roadside} and \cite{taymaz2020healthcare} seek to offer different service packages to patients; and \cite{almeida2021two-step} tackle the assignment of both the equipment and the additional hours of specialists to the located facilities. \cite{srivastava2021strengthening} use a maximal covering location problem to determine the optimal location of cold chain points and the allocation of vaccination session sites to these points. Moreover, \cite{de2022optimization}, \cite{decampos2024multi}, and \cite{defreitas2024integrated} analyze how to optimally locate MRI equipment, MS units, and NS technology at existing healthcare facilities, respectively, as well as the allocation of the demand to be served at those facilities.

Routing decisions play a crucial role in many optimization problems, and they are often integrated with location decisions to address real-world challenges effectively. This is the case in \cite{doerner2007multicriteria}, where the aim is to determine both the optimal location of a single mobile facility and the most efficient route for it to serve a set of demand nodes. Both p-median and maximal covering models are used to determine facility locations. \cite{vonachen2016optimizing} consider the location and assignments of community health workers (CHWs), as well as the location, routing, and assignments of the health workers' supervisors. For the first problem, the authors develop a capacitated set covering model (\citealp{daskin1997}; \citealp{guha2002}), while for the second one, the location-routing problem described by \cite{berger2007location} is utilized. Similarly, \cite{cherkesly2019community} define the location-routing covering problem, a variant of the location-routing and the covering tour problems, to design a community healthcare network. Decisions involve determining the location of CHWs and their supervisors, the allocation of CHWs to demand areas, and the supervisors' routes to visit and train CHWs. \cite{boutilier2020ambulance} adopt a two-stage approach to first strategically locate ambulance outposts and then optimize the routes of these emergency vehicles from the outposts to demand points. They integrate optimization techniques with machine learning approaches and simulation modeling. \cite{yang2021outreach} combine the set covering problem and vehicle routing problem to optimize the deployment of a team of clinicians tasked with setting up mobile clinics in hard-to-reach locations. While  \cite{alban2022resource} also consider mobile facilities, their model exclude routing decisions and instead focuses on the long-term allocation of mobile healthcare unit resources (visit capacities) to a large number of sites. The authors formulate this resource allocation problem as a knapsack problem with a sum of sigmoidal functions as the objective, each representing the demand for health service adoption at a site.

In specific contexts of developing countries, having the necessary funds to establish new medical facilities or relocate existing ones can be challenging. Some authors have focused on improving accessibility to existing facilities solely through enhancements in the transportation network (2 out of 56 papers). This approach is addressed by \cite{murawski2009improving}, who introduce the maximal covering network improvement model. In their work, they consider upgrading the transport links to ensure all roads are passable throughout the year. \cite{heyns2024optimisation} generalize the previous work by considering not only links for vehicular transportation but also walking, recognizing the importance of pedestrian mobility in areas where using motorized vehicles may not always be feasible.
In other studies (6 out of 56 papers), investments in road construction or improvements are sometimes combined with facility location (\citealp{cocking2012improving, ghaderi2013modeling, mohammadi2014design, shishebori2015robust, chen2022multi-objective, pourrezaie-khaligh2022fix-and-optimize}).

Overall, the reviewed literature shows a clear trend toward addressing long-term healthcare accessibility through location problems, often integrated with patient allocation decisions. While the allocation of equipment or services is also frequently considered, the explicit allocation of healthcare personnel to facilities receives comparatively less attention. Furthermore, integrated decision frameworks that jointly address location, allocation, and routing or network design reflect an increasing recognition of the interdependence between facility location and transport means or infrastructure.

\subsubsection{Mathematical modeling}\label{subsec:math_model}
The fourth column of Table~\ref{tab:opt_approach} shows that the majority of the selected papers (37 out of 56) use either integer linear programming (ILP) (e.g., \citealp{doerner2007multicriteria}; \citealp{yao2014locational}; \citealp{vonachen2016optimizing}; \citealp{mendoza-gomez2022location, mendoza-gomez2022regionalization,mendoza2024maximal,mendoza2024regionalization}) or mixed-integer linear programming (MILP) (e.g., \citeauthor{smith2009planning}, \citeyear{smith2009planning}, \citeyear{smith2013bicriteria}; \citealp{shariff2012location}; \citealp{gunecs2014}; \citealp{nunez2016column}; \citealp{chouksey2022hierarchical, chouksey2022optimisation,chouksey2024accelerated}) models. There are some works (13 out of 56) in which a mixed-integer nonlinear programming (MINLP) model is presented (\citealp{galvao2006load}; \citealp{ghaderi2013modeling}; \citealp{zahiri2014multi-objective}; \citealp{beheshtifar_multiobjective_2015}; \citealp{mousazadeh_accessible_2018,mousazadeh2018health}; \citealp{zarrinpoor2018design};  \citealp{de2020roadside}; 
 \citealp{motallebi2020short-term}; \citealp{pourrezaie-khaligh2022fix-and-optimize}; \citealp{rouhani2022robust}; \citealp{decampos2024multi}; and \citealp{defreitas2024integrated}). All of them, except for \cite{beheshtifar_multiobjective_2015} and \cite{decampos2024multi}, provide either a linear counterpart or a convex model for their MINLP formulations. In \cite{alban2022resource}, nonlinear programming (NLP) involving sigmoidal functions is used to model adoption dynamics. In particular, the authors adapt the Bass model \citep{bass1969} to capture these dynamics and estimate the adoption levels resulting from different allocation strategies.

To cope with non-deterministic features, different approaches have been adopted in the literature. One of them is fuzzy programming (FP) (5 out of 56 papers), which offers a solution for managing epistemic uncertainty in input parameters (possibilistic programming) and elasticity in constraints and/or flexibility in goals (flexible programming). \cite{zahiri2014multi-objective} 
 present a possibilistic programming model, and \cite{zahiri2014robust} and \cite{mousazadeh2018health} introduce different robust possibilistic chance-constraint programming approaches. \cite{mousazadeh_accessible_2018} combine possibilistic programming with flexible programming, providing a robust mixed possibilistic-flexible programming approach. Some works have adopted a stochastic programming (SP) approach (5 out of 56 papers)  (\citealp{mohammadi2014design}; \citealp{khodaparasti2018multi-period}; \citealp{wang2018healthcare}; \citealp{taymaz2020healthcare}; \citealp{yang2021outreach}). In particular, \cite{mohammadi2014design} combine FP with SP, introducing the mixed interval-fuzzy stochastic programming method. The proposed model can effectively manage mixed uncertainties represented as intervals as well as possibilistic and probabilistic distributions. \cite{khodaparasti2018multi-period} introduce a multi-period probabilistic stochastic programming model with chance constraints. These authors also provide an equivalent integer nonlinear programming (INLP) model with its corresponding linearization. 
\cite{wang2018healthcare} propose a bi-level programming (BLP) model that integrates stochastic and fuzzy factors. Meanwhile, \cite{andrade2015abc} present a model that integrates the probabilistic methodology introduced by \cite{marianov1996} with the dynamic model proposed by \cite{schmid2010} to create a unified stochastic dynamic programming (SDP) formulation of their location-allocation problem. 

Robust optimization (RO) is widely adopted as an effective approach for addressing uncertain optimization problems (6 out of 56 papers). For instance, \cite{shishebori2015robust} apply a scenario-based approach \citep{mulvey1995} to cope with uncertainties and a RO method with augmented robust constraints to deal with system's disruptions. \cite{zarrinpoor2017design} present a two-stage RO model in which the first-stage solution
is robust against potential disruption scenarios that are only revealed in the second stage, enabling demand reallocation when facilities are affected, and 
\cite{zarrinpoor2018design} formulate both a mixed-integer quadratically constrained quadratic programming model (with its corresponding convex model) and a robust scenario-based stochastic programming model. \cite{boutilier2020ambulance} provide a two-stage RO model along with an equivalent single-stage MILP, developing an edge-based reformulation of the classical path-based p-median problem. 
\cite{motallebi2020short-term} present a MINLP model that incorporates queuing theory to represent short-term uncertainty in demand and service times. They derive a linearized version and later propose a robust counterpart model to handle long-term demographic uncertainty. Meanwhile, \cite{rouhani2022robust} begin with a deterministic MINLP formulation, convert it into an equivalent linear form, and subsequently extend it through a RO approach to cope with uncertainty in demand and supply.

Finally, \cite{dogan_model_2020} present a goal programming (GP) approach to tackle a multi-objective problem, while \cite{khodaparasti2017enhancing} propose a fuzzy goal programming (FGP) method to address ambiguity in aspiration levels.

\subsubsection{Objective function}\label{subsec:obj_function}
To tackle conflicting objectives such as maximizing coverage or minimizing costs, researchers typically employ multi-objective optimization models. Column 5 of Table~\ref{tab:opt_approach} details how the objective function of the selected papers is formulated, indicating whether it is a single-objective (S) or a multi-objective (M). In this latter case, the number of objectives along with the method used to address them is specified, represented as M/number\_of\_objectives/method. 

Different techniques have been developed to solve multi-objective problems. Among these, scalar approaches such as the weighted sum (WS) method or the $\varepsilon$-constraint method are popular strategies to define an adequate single-objective problem, although they differ significantly in how they handle the multiple objectives. The WS method combines all objectives into a single composite function through a weighted sum, requiring the decision-maker to specify the relative importance of each objective. In contrast, the $\varepsilon$-constraint method selects one objective to optimize while transforming the remaining objectives into constraints with specified bounds, thereby allowing the generation of Pareto-optimal solutions. Among the selected papers, \cite{smith2013bicriteria}, \cite{nunez2016column}, \cite{de2020roadside}, \cite{karakaya2022biobjective}, \cite{eksioglu2024designing}, and \cite{heyns2024optimisation} use the classic weighted sum method. 
\cite{galvao2006load} apply the $\varepsilon$-constraint method introduced by \cite{cohon1978}.
\cite{pourrezaie-khaligh2022fix-and-optimize} and \cite{rouhani2022robust} employ the augmented $\varepsilon$-constraint method introduced by \cite{mavrotas2009}, while \cite{mousazadeh_accessible_2018} and \cite{mendoza-gomez2022location} utilize the improved version of \cite{mavrotas2013}.

Alternatively, certain heuristic methods maintain the objectives separately throughout the solving process, without aggregating them into a single-objective function. \cite{doerner2007multicriteria}, \cite{beheshtifar_multiobjective_2015}, \cite{wang2018healthcare}, \cite{motallebi2020short-term}, \cite{decampos2024multi}, and \cite{heyns2024optimisation} use evolutionary algorithms (EA) to address the multi-objective nature of their problems, aiming to approximate the set of Pareto-efficient solutions. In particular, \cite{doerner2007multicriteria} employ a generalization of the ant colony optimization metaheuristic \citep{dorigo1996} to treat the location and routing aspects of their problem simultaneously. These authors also apply two genetic-based algorithms that address these aspects separately, the vector evaluated genetic algorithm \citep{schaffer2014} and the multi-objective genetic algorithm  \citep{fonseca1993}, and perform comparisons between these three approaches. \cite{beheshtifar_multiobjective_2015}, \cite{motallebi2020short-term}, \cite{decampos2024multi}, and \cite{heyns2024optimisation} adapt the non-dominated sorting genetic algorithm II \citep{deb2002} to identify the Pareto-optimal set, and \cite{decampos2024multi} also apply the strength Pareto evolutionary algorithm \citep{zitzler2001}. To support decision-making, \cite{beheshtifar_multiobjective_2015} additionally use the \textit{a posteriori} preference method TOPSIS to rank these Pareto-optimal candidates and assist in selecting a most preferred compromise solution based on different weight vectors.
\cite{wang2018healthcare} incorporate four objective functions and develop a bi-level multi-objective particle swarm optimization method to simultaneously address location decisions and capacity adjustments. \cite{chen2022multi-objective} apply the improved multi-objective simulated annealing (SA) algorithm introduced in \cite{chen2021optimizing} for optimizing the location of rural EMS facilities.

Other strategies to cope with multi-objective problems are game theory (GT) \citep{mohammadi2014design}, goal programming (GP) \citep{dogan_model_2020}, or fuzzy approaches (\citealp{zahiri2014multi-objective}; \citealp{mousazadeh2018health}). The latter works use the \cite{Torabi2008}'s method (TH), which consists of an interactive possibilistic programming approach to find a preferred compromise solution. \cite{khodaparasti2017enhancing} provide a fuzzy goal programming (FGP) method that integrates the fuzzy set theory into the GP, allowing decision-makers to express target levels in an imprecise manner.

Most of the reviewed papers (34 out of 56) model their problems with a single objective function. Among those using multi-optimization approaches (22 out of 56), there is a clear preference for scalarization techniques and $\varepsilon$-constraint methods (11 out of 22 papers), probably because their conceptual simplicity and compatibility with exact optimization solvers. At the same time, there is growing use of evolutionary and metaheuristic approaches (7 out of 22) to preserve the multi-objective structure and approximate Pareto-efficient solutions in large-scale or highly complex settings.

\subsubsection{Solution approach}\label{subsec:sol_approach}
The last two columns in Table~\ref{tab:opt_approach} specify the solution approach(es) adopted by the selected papers. Many studies rely on commercial solvers (CS) like CPLEX (e.g., \citealp{murawski2009improving}; \citealp{shishebori2015robust}; \citealp{dogan_model_2020}; \citealp{mendoza2024maximal}), BARON (\citealp{mohammadi2014design}; \citealp{zahiri2014multi-objective}), or Gurobi (\citealp{kunkel2014optimal}; \citealp{yang2021outreach}; \citealp{chouksey2022hierarchical, chouksey2022optimisation,chouksey2024accelerated}), among others, to solve the proposed mathematical models. Only eight works utilize custom-built exact procedures: \cite{zarrinpoor2017design} and \cite{eksioglu2024designing} apply a Benders decomposition (BD) algorithm, while \cite{zarrinpoor2018design} and \cite{chouksey2024accelerated} develop an accelerated BD method; \cite{khodaparasti2017enhancing} apply the outer approximation (OA) algorithm to obtain the optimal nursing homes' sites; \cite{cherkesly2019community} and \cite{boutilier2020ambulance} develop variable generation (VG) and row and column generation (RCG) techniques, respectively, to reduce the dimensions of their problems and solve them efficiently; and \cite{flores2021optimizing} develop an ad-hoc brute force algorithm.

To solve large instances, some authors propose approximate approaches. As already mentioned in the previous section, some works tailor different genetic algorithms (GA) (\citealp{doerner2007multicriteria}; \citealp{beheshtifar_multiobjective_2015}; \citealp{motallebi2020short-term};  \citealp{decampos2024multi}; \citealp{heyns2024optimisation}), ant colony (AC) algorithms \citep{doerner2007multicriteria}, particle swarm optimization (PSO) algorithms \citep{wang2018healthcare}, other evolutionary algorithms (EAs) \citep{decampos2024multi}, or simulated annealing (SA) methods \citep{chen2022multi-objective} to solve their multi-objective models. \cite{shariff2012location} also propose a modified GA to solve the capacitated maximal covering problem. Other metaheuristic (MH) algorithms have been developed in several works. For instance, \cite{mendoza2024regionalization} propose an iterated greedy algorithm (IGA) to solve the model originally introduced in \cite{mendoza-gomez2022regionalization}. \cite{ghaderi2013modeling} and \cite{pourrezaie-khaligh2022fix-and-optimize} design fix-and-optimize (FO) strategies based on SA and tabu search, respectively, along with CPLEX. In addition to SA methods, \cite{mohammadi2014design} and \cite{zahiri2014multi-objective} also develop population-based evolutionary techniques such as an imperialist competitive algorithm (ICA) and a self-adaptive differential evolution (DE) algorithm, respectively. \cite{andrade2015abc} propose an artificial bee colony (ABC) algorithm to guide the mobile EMS of S\~{a}o Paulo in decisions regarding the location of ambulance stations and allocation and repositioning of these vehicles. 
Lagrangean (Lag) (\citealp{galvao2006load}; \citealp{alban2022resource}; \citealp{karakaya2022biobjective}), greedy (Gr) (\citealp{ghaderi2013modeling}; \citealp{srivastava2021strengthening}), column generation (CG) \citep{nunez2016column} and scenario generation (SG) \citep{boutilier2020ambulance} based, Banders-type \citep{chouksey2024accelerated}, and ad-hoc heuristics (H) (\citealp{cherkesly2019community}; \citealp{alban2022resource};  \citealp{kumar2023location}) have also been designed. In addition to a Lagrangean relaxation approach to solve their proposed integrated hierarchical location-allocation model, \cite{karakaya2022biobjective} also apply two top-down heuristics (TDH). Essentially, these approaches first locate upper-level facilities and then move down the hierarchy to locate lower-level facilities by sequentially solving different p-median models. 
Finally, \cite{chouksey2022hierarchical} present a sequential approach using relax-and-fix (RF) and FO heuristics, \cite{chouksey2024accelerated}  design a hybridized FO and SA strategy, and \cite{mohammadi2014design} develop a hybrid solution approach (HSA) using queuing theory, game theory, and interval, fuzzy, and stochastic programming techniques.

Decision processes are sometimes involved when dealing with multi-objective approaches and represent somehow the last tier of the solution phase. For instance, \cite{beheshtifar_multiobjective_2015} propose a technique for order of preference by similarity to ideal solutions (TOPSIS) for the decision-maker to find a preferred solution from the previously-obtained Pareto-efficient ones, based on their proximity to ideal reference points.


\subsection{Practical implementation} \label{subsec:pract_impl}
This section analyzes the real-world implementation of the proposed models or the absence thereof. The success of a study in terms of practical implementation is linked to 1) the use of real-world data as opposed to illustrative data, and 2) the involvement of stakeholders throughout the life cycle of the project, including problem definition, conceptual model building, validation and experimentation, and implementation (see \citealp{tako2015}; \citealp{brailsford2019}; \citealp{noorain2024}). For each case study, Table~\ref{tab:implementation} scrutinizes the input data sources, stages and methods of stakeholder or end-user engagement, and the level of real-world implementation.

\subsubsection{Model input sources} \label{subsec:mod_inputs}

To define the source for input data, we rely on the conceptual framework by \cite{brailsford2019}: 
\begin{enumerate}
    \item Real-world data (RD): data from case study organizations (including primary data collection and retrieval of data from information systems) or secondary data, in which some input data may also include estimates relying on context-specific features of the real-world setting.
    \item Illustrative data (ID): 
    data used purely for demonstration or model testing purposes with no justification. This includes synthetic, assumed, or randomly generated data, often based on expert estimates or hypothetical scenarios.
    \item Mixed data (MD): a combination of real-world data and illustrative data.
\end{enumerate}

\newgeometry{left=2cm, right=1cm, bottom=0.2cm, top=0.7cm}

\begin{table}[ht!]\scriptsize
\centering
\caption{Practical implementation aspects for the selected papers.}\label{tab:implementation}
\resizebox{0.95\linewidth}{!}{
\hspace{-1.2cm}\begin{tabular}{p{0.5cm}p{5.8cm}!{\color{gray!30}\vrule}p{1.6cm}!{\color{gray!30}\vrule}p{2.5cm}!{\color{gray!30}\vrule}p{2cm}!{\color{gray!30}\vrule}p{1cm}!}
\noalign{\smallskip}\hline\noalign{\smallskip} 
No. & Reference & \parbox{3cm}{Input data} & \parbox{4cm}{Stakeholder \\engagement \\stage} & \parbox{4cm}{Stakeholder \\engagement \\method} & \parbox{1.2cm}{Level of \\implementation} \\
\noalign{\smallskip}\hline\noalign{\smallskip}

    1 & \cite{galvao2006load}     &  MD   & PD & M   & PC    
    \\[0.2em]
     \arrayrulecolor{gray!30}\hline

    2 & \cite{doerner2007multicriteria}     &  RD    & PD  & U   & PC    
    \\[0.2em]
     \arrayrulecolor{gray!30}\hline
     
    3 & \cite{murawski2009improving}     &  MD   & N   &  -  &  PC   \\[0.2em]
     \arrayrulecolor{gray!30}\hline
     
    4 & \cite{smith2009planning}     & RD  & PD, MB    & M, Q  &  P \\[0.2em] \arrayrulecolor{gray!30}\hline
    
    5 & \cite{cocking2012improving}       &  RD   & PD  & I  & P \\[0.2em] \arrayrulecolor{gray!30}\hline
    
    6 & \cite{shariff2012location}       & MD    & N   & -    & PC   \\[0.2em] \arrayrulecolor{gray!30}\hline
    
    7 & \cite{ghaderi2013modeling}     & MD    &  N & -   & PC 
    \\[0.2em] \arrayrulecolor{gray!30}\hline

    8 & \cite{smith2013bicriteria}     &  MD    & PD  & U  & PC   
    \\[0.2em]
     \arrayrulecolor{gray!30}\hline
     
    9  & \cite{gunecs2014}     & MD   & PD  & M   & PC 
    \\[0.2em] \arrayrulecolor{gray!30}\hline

    10 & \cite{kunkel2014optimal}   &  RD   &  PD, VE &  M, I   & P  \\[0.2em] \arrayrulecolor{gray!30}\hline

    11 & \cite{mohammadi2014design}      & ID    & PD   & U     & PC \\[0.2em] \arrayrulecolor{gray!30}\hline

    12 & \cite{yao2014locational}      &  RD  & PD  & Q  & PC \\[0.2em] \arrayrulecolor{gray!30}\hline
    
    13 & \cite{zahiri2014multi-objective}      & ID    & PD & U  & PC \\[0.2em] \arrayrulecolor{gray!30}\hline

    14 & \cite{zahiri2014robust}      & ID    & PD  & U    &  PC  \\[0.2em] \arrayrulecolor{gray!30}\hline
    
    15 & \cite{andrade2015abc}     &  RD   & PD, VE, RI   & C     & E \\[0.2em] \arrayrulecolor{gray!30}\hline    
    
    16 & \cite{beheshtifar_multiobjective_2015}    &  MD & VE  &  R   & PC  \\[0.2em] \arrayrulecolor{gray!30}\hline

    17 & \cite{shishebori2015robust}      &   MD  & N   & -   & PC  \\[0.2em] \arrayrulecolor{gray!30}\hline

    18 & \cite{lim2016coverage}    &  MD & PD   & U   & PC  \\[0.2em] \arrayrulecolor{gray!30}\hline

    19 & \cite{nunez2016column}    &   ID, RD  & PD  &  I    &  PC \\[0.2em] \arrayrulecolor{gray!30}\hline    
    
    20 & \cite{vonachen2016optimizing}    &  RD   &  PD, MB, VE, RI & C     & E \\[0.2em] \arrayrulecolor{gray!30}\hline
    
    21 & \cite{khodaparasti2017enhancing}  &   RD  & PD, MB, VE    & C & P \\[0.2em] \arrayrulecolor{gray!30}\hline

    22 & \cite{zarrinpoor2017design}     &  MD   & PD  & U   & PC 
    \\[0.2em]
     \arrayrulecolor{gray!30}\hline
     
    23 & \cite{khodaparasti2018multi-period}     &   MD  & N  & -  & PC \\[0.2em] \arrayrulecolor{gray!30}\hline
    
    24 & \cite{mousazadeh_accessible_2018}    &  ID   & PD, VE
    & U, R
    &  PC\\[0.2em] \arrayrulecolor{gray!30}\hline

    25 & \cite{mousazadeh2018health}     &  ID    & PD  & C   & PC   
    \\[0.2em]
     \arrayrulecolor{gray!30}\hline
     
    26 & \cite{wang2018healthcare}    & RD &   PD   & U   & PC \\[0.2em] \arrayrulecolor{gray!30}\hline
    
    27 & \cite{zarrinpoor2018design}    &  MD   & PD     & U   & PC\\[0.2em] \arrayrulecolor{gray!30}\hline
    
    28 & \cite{cherkesly2019community}   & RD    & PD, MB, VE, RI  & C   & E \\[0.2em] \arrayrulecolor{gray!30}\hline
    
    29 & \cite{boutilier2020ambulance}   &  MD   & PD   & M, Q & P \\[0.2em] \arrayrulecolor{gray!30}\hline

    30 & \cite{de2020roadside}    & RD  &  PD, MB, VE, RI
    & C, I &  E  \\[0.2em] \arrayrulecolor{gray!30}\hline

    31 & \cite{dogan_model_2020}    &  MD &  PD, VE
    &  U, R
    &  PC \\[0.2em] \arrayrulecolor{gray!30}\hline    

    32 & \cite{motallebi2020short-term}    &  RD   & PD  &   M  & PC\\[0.2em] \arrayrulecolor{gray!30}\hline

    33 & \cite{taymaz2020healthcare}     &  MD   & PD, MB   &  C  & PC \\[0.2em] \arrayrulecolor{gray!30}\hline
    
    34 & \cite{almeida2021two-step}    & RD  & PD   & U   & PC  \\[0.2em] \arrayrulecolor{gray!30}\hline
    
    35 & \cite{flores2021optimizing}    &   RD    &  N  & -    & PC  \\[0.2em] \arrayrulecolor{gray!30}\hline
    
    36 & \cite{srivastava2021strengthening}    &  RD   & PD, MB, VE   &  C &  P \\[0.2em] \arrayrulecolor{gray!30}\hline
     
    37 & \cite{yang2021outreach}    &  ID & PD   & U    & PC \\[0.2em] \arrayrulecolor{gray!30}\hline

    38 & \cite{alban2022resource}     &  RD   & PD, MB  & C   & P    
    \\[0.2em]
     \arrayrulecolor{gray!30}\hline
     
    39 & \cite{chen2022multi-objective}   &  ID   &  N  &   -   &  PC \\[0.2em] \arrayrulecolor{gray!30}\hline

    40 & \cite{chouksey2022hierarchical}   &  MD   & N &   - & PC \\[0.2em] \arrayrulecolor{gray!30}\hline

    41 & \cite{chouksey2022optimisation}    &   MD  &  N &  -  &  PC \\[0.2em] \arrayrulecolor{gray!30}\hline
        
    42 & \cite{de2022optimization}    & MD   & N  &  -   &  PC\\[0.2em] \arrayrulecolor{gray!30}\hline

    43 & \cite{elorza2022assessing}    &   RD    &  PD, MB &  C     & P  \\[0.2em] \arrayrulecolor{gray!30}\hline

    44 & \cite{karakaya2022biobjective}     &  MD    & N  & -   & PC 
    \\[0.2em]
     \arrayrulecolor{gray!30}\hline
     
    45 & \cite{mendoza-gomez2022location}      &  RD   &  N  & - & PC \\[0.2em] \arrayrulecolor{gray!30}\hline

    46 & \cite{mendoza-gomez2022regionalization}   & RD  & N  &  - &  PC  \\[0.2em] \arrayrulecolor{gray!30}\hline

    47 & \cite{pourrezaie-khaligh2022fix-and-optimize} &  ID  &  N &  - &   PC      \\[0.2em] \arrayrulecolor{gray!30}\hline

    48 & \cite{rouhani2022robust}       &  RD   & PD  & U & PC \\[0.2em] \arrayrulecolor{gray!30}\hline
    
    49 & \cite{kumar2023location}     &  MD   &  N  & -    & PC \\[0.2em] \arrayrulecolor{gray!30}\hline

    50 & \cite{decampos2024multi}     &  RD    & N & -   & PC    
    \\[0.2em]
     \arrayrulecolor{gray!30}\hline

    51 & \cite{defreitas2024integrated}     &  MD    & N  & -   & PC   
    \\[0.2em]
     \arrayrulecolor{gray!30}\hline
     
    52 & \cite{eksioglu2024designing}   &   RD  & PD &  M &  PC \\[0.2em] \arrayrulecolor{gray!30}\hline

    53 & \cite{heyns2024optimisation}    & MD   &  VE &  M &  PC  \\[0.2em]
     \arrayrulecolor{gray!30}\hline

    54 & \cite{mendoza2024maximal}     &  MD    & N  & -  & PC   
    \\[0.2em]
     \arrayrulecolor{gray!30}\hline

    55 & \cite{mendoza2024regionalization}     &  RD   & N  & -   & PC    
    \\[0.2em]
     \arrayrulecolor{gray!30}\hline

    56 & \cite{chouksey2024accelerated}     &  MD   & N  & -   & PC   
    \\
        
     \arrayrulecolor{black}
\noalign{\smallskip}\hline\noalign{\smallskip}
    \multicolumn{6}{l}{\parbox{1\linewidth}{\textbf{Input data}: MD (mixed data), RD (real-world data), ID (illustrative data); \textbf{Stakeholder engagement stage}: PD (problem definition), N (none), MB (model building), VE (validation and experimentation), RI (real implementation); 
    \textbf{Stakeholder engagement method}: M (meetings, consultations, and other forms of communication), U (unclear), 
    Q (questionnaires, surveys), 
    I (interviews),
    C (close collaboration),  
    R (ranking);
    \textbf{Level of implementation}: PC (proof of concept), P (potential), E (evidence of real implementation). 
    }}
    \end{tabular}}
\end{table}

\restoregeometry


The third column of Table~\ref{tab:implementation} shows that just eight studies rely only on illustrative data sources. Interestingly, all but \cite{yang2021outreach} and \cite{chen2022multi-objective} have their case studies in Iran. In particular, \cite{zahiri2014multi-objective} and \cite{zahiri2014robust} rely on experts' opinions to locate candidate sites for health centers, and the authors of the former study also use estimated values for all parameters in their model. For instance, they consider the saving cost of integrated facilities based on the construction costs, which at the same time are estimated based on historical data. 
The latter study proposes several test problems to validate their model, the same as in \cite{yang2021outreach} with data adapted from four countries in sub-Saharan Africa. In \cite{mousazadeh2018health}, the authors also use estimations for all model parameters, most of them based on decision-makers' suggestions.

The majority of the papers use real-world data or mixed data sources for their case studies. There are several works employing commercial GIS software 
to record and process spatial data. For instance, \cite{doerner2007multicriteria} gather GIS information based on detailed regional maps and historical census data; \cite{murawski2009improving} extract road network representations and features, clinic locations, and village populations' information from geo-databases and use GIS functions to map the results;  \cite{yao2014locational} rely on ArcGIS to conduct spatial analysis and visual evaluations of their collected data, obtained through census records and surveys to patients; and \cite{vonachen2016optimizing} use GIS software to process GPS data on communities and transportation networks, which were manually collected by their partners. These data are later on processed using ArcGIS, and the authors emphasize the importance of using actual road network distances rather than Euclidean estimations (as is common practice).  Other studies using commercial GIS software are \cite{kunkel2014optimal}, \cite{beheshtifar_multiobjective_2015}, \cite{cherkesly2019community}, and \cite{chen2022multi-objective}. Open source geographical tools, such as Google Maps, Google Earth, OpenStreetMap, Mapbox, or OpenRouteService are also used to estimate travel times  (\citealp{zahiri2014robust}; \citealp{andrade2015abc}; \citealp{flores2021optimizing}; \citealp{srivastava2021strengthening}; \citealp{rouhani2022robust}; \citealp{defreitas2024integrated}), to map road networks (\citealp{vonachen2016optimizing}; \citealp{flores2021optimizing}; \citealp{chen2022multi-objective}), or to identify the coordinates of health facilities or villages (\citealp{flores2021optimizing}; \citealp{chouksey2022hierarchical,chouksey2024accelerated}; \citealp{defreitas2024integrated}). \cite{alban2022resource} obtain population data using the WorldPop project.

Only a few studies have directly collected data in the field (\citealp{cocking2012improving}; \citealp{vonachen2016optimizing}; \citealp{cherkesly2019community}; \citealp{boutilier2020ambulance}; \citealp{srivastava2021strengthening}), such as demographic information or transportation network characteristics. In some cases, data are provided first-hand by different organizations. For instance, \cite{smith2009planning} collect data from a hospital community health project in India; \cite{gunecs2014} obtain information from the Sakarya Provincial Health Department; \cite{kunkel2014optimal} use two datasets from different institutions originated from the Malawi Housing and Population Census; \cite{andrade2015abc} obtain historical data on call locations, times, and severity from a computer-aided dispatch system owned by the EMS in S\~ao Paulo, Brazil; \cite{de2020roadside} and \cite{taymaz2020healthcare} use data sources provided by North Star Alliance, a nongovernmental organization (NGO) operating in Africa; and \cite{alban2022resource} use data and documentation from MSI Reproductive Choices, an NGO that provides family planning services in Uganda and other countries.

Several works rely on secondary data from official databases. To name a few, \cite{shariff2012location},  \cite{dogan_model_2020} and \cite{karakaya2022biobjective}, and \cite{mendoza2024maximal} obtain population data from the Housing Census Malaysia, Turkish Statistical Institute, and Mexico Population and Housing Census, respectively;
\cite{de2022optimization} and \cite{decampos2024multi} consider the National Health System's database of Brazil to gather information on service capacity or demands; and \cite{rouhani2022robust} 
obtain data about the supply and demand of organs from the Transplantation and Special Disease Department of Iran. 

The lack of data in developing countries and the prevalence of significantly outdated databases makes it natural to consider estimated parameters. Among the most commonly estimated parameters, it is usual to take the coverage distance from the literature \citep{murawski2009improving} or directly specify it without any justification for experimental purposes 
(\citealp{smith2013bicriteria}; \citealp{shariff2012location}; \citealp{mendoza-gomez2022location}; \citealp{kumar2023location}; \citealp{mendoza2024maximal}). Nevertheless, some works set it based on goals established by healthcare institutions, policymakers, or expert suggestions (\citealp{doerner2007multicriteria}; \citealp{cocking2012improving}; \citealp{lim2016coverage}; \citealp{mousazadeh_accessible_2018, mousazadeh2018health}; \citealp{chouksey2022hierarchical, chouksey2022optimisation,chouksey2024accelerated}; \citealp{decampos2024multi}; \citealp{eksioglu2024designing}). The costs of opening new facilities are typically estimated based on historical data (\citealp{zahiri2014multi-objective}; \citealp{motallebi2020short-term}; \citealp{mendoza-gomez2022location}; \citealp{rouhani2022robust}) or according to the authors' best judgment (\citealp{shishebori2015robust}; \citealp{mousazadeh_accessible_2018}; \citealp{pourrezaie-khaligh2022fix-and-optimize}; \citealp{chouksey2022hierarchical, chouksey2022optimisation};  \citealp{kumar2023location}). This latter approach is usually applied to estimate investment or construction costs for road links  (\citealp{murawski2009improving}; \citealp{cocking2012improving}; \citealp{ghaderi2013modeling}; \citealp{shishebori2015robust}; \citealp{chen2022multi-objective}), the average speed of vehicles (\citealp{doerner2007multicriteria}; \citealp{chen2022multi-objective}; \citealp{pourrezaie-khaligh2022fix-and-optimize}), or the acquisition and maintenance of equipment \citep{defreitas2024integrated}. There are also occasions where some of these parameters are estimated based on information provided by local public health officials, as in \cite{gunecs2014}, \cite{mousazadeh2018health}, or 
\cite{yang2021outreach}.

\subsubsection{Stakeholder engagement} \label{subsec:stakehol_eng}
To understand the level of real stakeholders, policymakers, or end-users engagement, we investigate: 
\begin{enumerate}
    \item When they are involved. This is based on the different study stages (see \citealp{brailsford2019}):
    \begin{enumerate}
        \item None (N), if they are not involved at all. 
        \item Problem definition (PD), which includes data collection and estimating parameters, as well as offering insights into the relevant application domain. 
        \item Model building (MB), which involves the identification of objectives and/or constraints.
        \item Validation and experimentation (VE), which involves confirming whether the model is reasonable or practical, as well as testing and discussing various scenarios to assess sensitivity. 
        \item Real implementation (RI), including a discussion of all the findings, barriers, or recommendations that can or cannot be implemented. 
    \end{enumerate}
    \item How they are involved. Sometimes it is unclear (U) from the studies, but they might have been involved through:
    \begin{enumerate}
        \item Meetings, consultations, and other forms of communication (M). 
        \item Questionnaires or surveys (Q). 
        \item Interviews (I). 
        \item Close collaborations (C) with certain organizations, such as NGOs or health authorities. 
    \end{enumerate}
\end{enumerate}
Columns 4 and 5 of Table~\ref{tab:implementation} summarize this information for the selected papers. There are 20 out of 56 works with no engagement at all, and from the remaining ones, 21 studies have involved stakeholders only at the initiation of the project. The majority of them do not clearly specify how this engagement is performed. For instance, \cite{doerner2007multicriteria} set an upper bound for an acceptable walking distance as defined by a political decision-maker, and \cite{smith2013bicriteria}, \cite{lim2016coverage}, and \cite{yang2021outreach} are given certain data to inform their case studies. Other works rely on experts' opinions to consider uncertainty \citep{mohammadi2014design}, or to determine specific parameters, such as the number and candidate locations for new health centers (\citealp{zahiri2014multi-objective,zahiri2014robust}; \citeauthor{zarrinpoor2018design}, \citeyear{zarrinpoor2017design}, \citeyear{zarrinpoor2018design}; \citealp{almeida2021two-step}; \citealp{rouhani2022robust}), the minimum and maximum coverage radius of different facilities \citep{mousazadeh_accessible_2018}, or the service and waiting times \citep{dogan_model_2020}. When integrating multi-criteria decision-making approaches, experts' opinions are sometimes sought to define the criteria or to give them weights (\citealp{wang2018healthcare}).  
Additionally, \cite{mousazadeh_accessible_2018} and \cite{dogan_model_2020} 
also incorporate decision-makers' ranking of conflicting objectives during the validation and experimentation phase, a practice also identified in \cite{beheshtifar_multiobjective_2015} as the sole involvement of stakeholders.

Other studies document stakeholder engagement more explicitly. For example, \cite{galvao2006load} incorporate capacity constraints informed by input from municipality health officials regarding shortages in the system; \cite{gunecs2014} set possible ranges for model parameters in discussion with city health authorities;  \cite{motallebi2020short-term} determine a minimum service level from conversations with high-level policymakers; \cite{mousazadeh2018health} report effective and close collaboration with the national and provincial healthcare experts from Iran's Ministry of Health and Medical Education
(MoHME), who provided the estimations or insights used to estimate model parameters; 
\cite{eksioglu2024designing} seek guidance from drone experts to obtain estimated costs and capabilities. 
In their project, \cite{cocking2012improving} conduct interviews early on to understand how road quality impacts ease of travel and its consequent effects on residents' visits to health centers. \cite{nunez2016column} rely on interviews with staff members of the NGO North Star Alliance to inform their problem. \cite{yao2014locational} survey women to collect geographic, demographic, and health factors that could inform the establishment of HIV testing facilities in rural Mozambique. Furthermore, \cite{boutilier2020ambulance} conduct surveys to potential patients in Dhaka, Bangladesh, to investigate their views on ambulance response times, and also consult with a transportation engineer to precisely identify the roads suitable for these vehicles. 

Besides including stakeholders and potential users in the initial phase of the case study, \cite{smith2009planning}, \cite{taymaz2020healthcare}, \cite{alban2022resource}, and \cite{elorza2022assessing} also engage them in shaping the design of their models. Specifically, \cite{smith2009planning} start by collecting patient data through house-to-house surveys and then define different equity and efficiency objectives for their location models according to contacts and discussions with medical personnel and members of faith-based organizations involved in community health initiatives. \cite{taymaz2020healthcare} assign communicable diseases to the corresponding service packages based on the interventions required, 
following consultations with staff from North Star Alliance, who also contribute most of the relevant data. This collaboration also yield precise definitions of coverage for different diseases to ensure compliance with continuum of care standards and aided in establishing certain parameters related to risk-averse measures within the proposed models. 
\cite{alban2022resource} weight adopter demand at each site and define the total capacity as well as lower and upper-bound constraints based on discussions with the NGO they collaborated with.  
Furthermore, \cite{elorza2022assessing} forge a collaboration with health authorities in the Bahia Blanca district of Argentina, who supply the authors with valuable data on healthcare services, system, and demographics, and offer insights on potential variables or parameters to inform their modeling. Meanwhile, \cite{kunkel2014optimal} receive data from various academic institutions and set additional parameters based on interviews with health surveillance assistants. Moreover, the results of their optimization recommendations are validated by stakeholders of the community health scheme under investigation. They identify certain travel distances as impractical, which influence their decision to designate additional facilities as backpack resupply centers when planning the program in the field. Subsequent to designing their model and conducting several experiments, \cite{heyns2024optimisation} seek expert feedback from economists, geographers, and transportation specialists through a series of personal meetings. These professionals value the optimization-driven outcomes compared to the existing road network and offer recommendations for future considerations.

\cite{khodaparasti2017enhancing} engage stakeholders during the first three stages of their case study through close collaboration with a community-based organization. They conduct field-based research together with health practitioners to obtain data, and incorporate the ambiguity of stakeholders' preferences and economic concerns into their model. Several scenarios are discussed to address research questions pertinent to decision-makers, which help them evaluate system improvements proposed by their approach compared to the current situation. On the other hand, in the research conducted by \cite{srivastava2021strengthening}, all authors are policymakers or stakeholders who directly contribute to the study. They acquire facilities' coordinates using a combination of tele-calling and site visits. Relevant information for the case study is gathered by volunteers and healthcare staff in the field, with subsequent validation and correction by project managers.

The only studies to expand stakeholders' engagement into the implementation phase are the ones by \cite{andrade2015abc}, \cite{vonachen2016optimizing}, \cite{cherkesly2019community}, and \cite{de2020roadside}. These four works forge distinct partnerships that remain actively involved across most stages of their projects. Specifically, \cite{andrade2015abc} collaborate with S\~ao Paulo's Mobile Emergency Attendance Service (SAMU-SP), which provides data and field experience to adjust several parameters in the model. Furthermore, SAMU-SP officials are engaged in the decision process to define the systems' expansion for the next years, which ends up impacting on the existing strategic, tactical, and operational plans on the allocation and relocation of their emergency vehicles. \cite{vonachen2016optimizing} work closely with the NGO Last Mile Health to help them improve their operational decisions regarding a community health initiative in rural Liberia. The authors develop a model based on the available manually collected data and contextual restrictions provided by the NGO. Furthermore, their partners contribute feasible initial solutions to assist in solving the model with commercial software, and ongoing joint efforts were being made for implementation at the time of writing their paper. These efforts culminate in \cite{cherkesly2019community}, where a decision-making tool is developed for the NGO to evaluate the costs and service level impacts of different varying parameters. The findings facilitate the development of the community healthcare network design in the field, which helps to identify barriers encountered during the implementation. Finally, \cite{de2020roadside} collaborate with the NGO North Star Alliance to gather data and gain a practical understanding of the problem under study. By means of interviews with staff members, the authors define the model objectives and constraints, which sometimes depend on NGO donors' preferences. A software package containing the model is developed jointly with the Dutch consulting group ORTEC 
and then installed at the headquarters.

\subsubsection{Level of implementation} \label{subsec:lev_impl}
This section attempts to identify whether the findings of the case studies have been implemented in practice or not. The model presented in each selected paper is categorized according to the following classification (see \citealp{brailsford2019}):
\begin{enumerate}
    \item Proof of concept (PC): the model is a proof of concept only.
    \item Potential (P): the model has been developed with stakeholders.
    \item Evidence of real implementation (E): the model has been developed and results implemented. This has influenced real-world decision-making.
\end{enumerate}

The last column of Table~\ref{tab:implementation} summarizes this information. Around 79\% of the models are simply proof of concept even if they had stakeholders' or end-users' engagement in the early stages of their case studies. There are only eight papers where the potential value of their optimization strategies is evidenced, as stakeholders were more actively involved throughout the project's evolution. For instance, \cite{smith2009planning} highlight that their contacts with health professionals give credence to their findings and perform a comprehensive sensitivity analysis that can provide local decision-makers with a range of options when setting up their community health program.  
Similarly, \cite{alban2022resource} emphasize the potential benefits of accounting for the future evolution of adopter demand when allocating resources, and demonstrate how such demand can be predicted without additional data-gathering efforts. They also offer managerial insights for MSI Reproductive Choices on how to allocate their resources based on predicted adoption dynamics, recommendations that sometimes diverge from the NGO's current intuitive practices. 
Some of these works also present reasons that make it challenging to bring their models into practice. In particular, \cite{cocking2012improving} emphasize that the primary obstacle to implementing their location and network design problem is the administrative hurdles involved, as the policymakers responsible for decisions on health facilities may differ from those overseeing decisions on road infrastructure. \cite{kunkel2014optimal} also identify political willpower as a factor that might contribute to deviations between current decisions and those suggested by their model, and stress the importance for their partners to consider these additional issues when putting their project into action. 
In \cite{boutilier2020ambulance}, researchers conduct fieldwork and directly observe situations occurring during ambulance responses. This facilitates their exploration of several scenarios that could assist ambulance providers or even the Bangladesh government in enhancing their systems. The authors conclude that a complete redesign might be necessary, but they also emphasize the requirement for additional initiatives to enable such changes, including awareness about emergency medical transport and central calls. 
\cite{srivastava2021strengthening} point out the lack of operational knowledge at the local level as a hindrance to the implementation of their model in practice, despite remote validation by project managers. The researchers highlight the need to engage with local medical staff to ensure that the results are accepted and can be effectively implemented, thereby positively impacting the last tier of the Indian immunization supply chain. 
\cite{elorza2022assessing} are confident that their system redesign proposal can potentially be implemented, but they acknowledge the difficulty of accurately estimating certain model parameters, even after consulting with specialists from various medical services.

As indicated in Section~\ref{subsec:stakehol_eng}, stakeholder engagement extends to the final implementation phase only in four papers, which are the sole works demonstrating varying degrees of real impact evidence. In the case of \cite{de2020roadside}, the developed models are integrated into software installed for the NGO staff overseeing the project. Nonetheless, while the paper mentions stakeholders' ability to adjust certain parameters within this tool package, it does not furnish additional evidence indicating implementation in the field. Similarly, \cite{vonachen2016optimizing} underscore the direct impact their models have on practice, as at the time of publication the NGO Last Mile Health had just initiated the implementation efforts of the community health program under consideration according to their results. However, it is difficult to ascertain the true extent of their findings on the ground, as these outcomes are unreported in the paper. 
Drawing from the previous research, \cite{cherkesly2019community} show some of the outcomes regarding the implementation of the community health plan. The project was in the design phase at the time of the paper writing, and the reported results pertain to the initial process of recruiting, locating, and routing both the medical personnel visiting remote communities and their supervisors. The project partners can adjust various parameters of the model, enabling them to address new arising situations, such as the potential unavailability of a worker, thereby having a temporally adaptive methodology. While the model indeed proposes alternative solutions for when a worker cannot be identified in a specific community, the NGO staff faces challenges in implementing the modifications suggested by the authors, associated with the difficulty of communicating with remote regions of Liberia. 
Perhaps the most well-documented work concerning the real impact of its practical implementation is \cite{andrade2015abc}. Despite initial reluctance and skepticism from authorities and the general society, SAMU-SP officials embraced an expansion plan based on the research findings, as these clearly demonstrated the benefits it entailed. Specifically, the decision was made to increase the number of mobile stations and to revamp the current operational plan for ambulance allocation and relocation according to the proposed mathematical model and heuristic. Ambulance response times and associated costs were substantially reduced while maintaining the same level of service. The positive outcomes of these newly introduced changes resulted in SAMU-SP receiving a distinction and the research gaining widespread public attention after being published in a prestigious Brazilian newspaper. At the time of article publication, SAMU-SP was actively utilizing the decision-making tool developed by the authors to further enhance their emergency vehicle system.

\section{Roadmap for future research}\label{sec:roadmap}

Optimizing healthcare accessibility in developing countries requires addressing the unique challenges these countries face, such as limited resources, geographic disparities, and socio-economic inequities. Analyzing the 56 selected papers, we have identified various gaps and opportunities for future research that OR practitioners should consider. These include new features to be modeled, optimization methods, and tools that facilitate the practical implementation of the proposed models.

Importantly, the gaps identified in this roadmap are not only derived from what has been addressed (or neglected) in the existing literature, but also from persistent challenges observed in practice. Thus, several of the research directions discussed below are additionally motivated by recurring real-world problems reported by international organizations or from insights from our collaborations with NGOs operating in developing countries.

\subsection{Novel features to be modeled}

Based on both the frequency with which challenges are reported in the analyzed papers and feedback from our NGO partners involved in healthcare delivery in Africa, the features discussed below are presented in approximate order of practical urgency. In particular, workforce-related decisions and service integration are consistently identified by the NGO staff as key concerns. Hence, this prioritization reflects applied needs rather than purely methodological novelty.

\subsubsection{Specialized services and staff allocation}

Strengthening service delivery and ensuring access to specialized care are central components of UHC. Recent WHO reports emphasize that progress towards UHC remains highly uneven across service domains, with particularly large gaps in service capacity, access to non-communicable disease (NCD) treatment, and the availability and distribution of skilled health workers in low- and middle-income countries \citep{WHO2025}. 
Persistent shortages of trained personnel, skill mismatches, and unequal geographical distribution of healthcare workers continue to constrain the effective delivery of specialized services, especially in regions where health needs are greatest. These challenges highlight the practical need for decision-support tools that jointly address the provision of multiple healthcare services and the allocation of specialized staff within health systems.

Our analysis of the 56 papers shows that the vast majority of the studies focus on a single level of medical care, with general primary healthcare or its various sub-disciplines being the most prominent (see Section~\ref{subsec:problem_features}). Secondary care is significantly underrepresented in the literature, despite its critical role in addressing NCD treatments and other conditions that account for a substantial share of the global disease burden, indicating ample opportunities for research on enhancing access to specialized medical services. Furthermore, as we previously pointed out, unlike in developed countries, many developing countries allow direct access to specialists without requiring a referral, making it even more natural to study this level of care independently. Although some works have focused on specific medium complexity procedures (e.g., \citealp{de2022optimization}) or services \citep{mendoza2024maximal} offered at the secondary care level, only \cite{almeida2021two-step} study a secondary healthcare system as a whole. This means that there are different specialties and, consequently, different types of health professionals. 
Their model informs about potential additional hours that specialists need to be hired in certain regions to meet existing demand. Nevertheless, there are several open issues to address, such as explicitly modeling the preferences and allocation of professionals to health centers. 
To the best of our knowledge, \cite{gunecs2014} is the only selected paper that integrates location decisions to improve patients' healthcare accessibility with staff allocation, although it is limited to family doctors. As mentioned above, workforce shortages are a growing problem globally, but the lack of specialized health workers is of even greater concern in low- and middle-income countries \citep{Liu2017}. Models that incorporate the location of multi-service health facilities along with decisions on both patient allocation and the allocation of distinct medical specialists are needed. 

Additionally, workforce allocation has been completely overlooked in papers proposing hierarchical models. Although these models address the diverse nature of care provided at various levels of health systems, they fail to consider the required skills of healthcare personnel and the associated costs. There is therefore ample opportunity for operational researchers to integrate tactical and operational decisions concerning medical staffing with strategic location decisions at different levels. 
For example, in the Gambella region of Ethiopia, the NGO Doctors with Africa CUAMM is interested in improving the healthcare system, which consists of various types of facilities (health posts, health centers, and hospitals) offering a range of services. In addition to determining where to open new facilities, one of their main challenges is the allocation and scheduling of family doctors, nurses, midwives, and other specialists to each facility, considering their preferences regarding distances to the facilities, as well as linguistic and cultural barriers. Their strategic decisions are inevitably constrained by the available health workers, making it essential to account for these factors when developing optimization models \citep{cuamm2023}.

Furthermore, it is important to highlight that many healthcare systems in these underdeveloped contexts are funded by donations from NGOs or other organizations. Stakeholders may decide to invest in training new skilled professionals if provided with an adequate and realistic decision-making tool.

\subsubsection{Integration of different service delivery channels}
Regarding the channels by which healthcare services are delivered to patients, we have distinguished between fixed facilities, mobile facilities, and health workers who travel to provide medical care, as discussed in Section~\ref{subsec:service_delivery_channel}. Frequently, studies on mobile clinics completely ignore preexisting fixed facilities (e.g., \citealp{doerner2007multicriteria}), thus oversimplifying reality and preventing these models from being applied effectively on the ground. Healthcare provision in underserved settings, however, often relies on a combination of fixed and mobile facilities. For example, the French Red Cross has, since 1987, operated a volunteer-based primary care model that integrates fixed access points with mobile structures to ensure access to healthcare for people experiencing severe social exclusion \citep{redcross2022}. In developing contexts, it is of particular importance to consider situations that integrate both types of service delivery channels, especially in scenarios where complete population coverage is required (such as routine vaccination). In such cases, two-stage models could be developed: the first phase would determine the locations of certain fixed facilities, while the second phase would address the remaining uncovered demand by establishing the stops and routing of mobile clinics.

Beyond physical service delivery, telemedicine has emerged as an additional channel of providing healthcare, particularly in response to pandemics such as COVID-19. Although its implementation remains controversial due to differences in medical consumption patterns and its overall impact on hospital operations \citep{zhou2023}, telemedicine projects have been successfully implemented in developing countries \citep{Combi2016}. These projects can help address medical shortages by providing remote consultations and extending healthcare access, especially in areas with limited access to specialists. However, several barriers hinder its widespread adoption, including poor internet infrastructure, unreliable electricity, high costs of technology, and limited digital literacy among both patients and healthcare providers. Despite these challenges, telemedicine has the potential to reduce pressure on overburdened facilities and reshape healthcare planning. This highlights the need to develop hybrid optimization models that integrate telemedicine with physical infrastructure, ensuring not only the strategic placement of facilities but also the efficient allocation of medical staff and resources to balance in-person and remote services effectively.

\subsubsection{Comprehensive multi-period settings}
Although some of the analyzed papers incorporate time as a modeling dimension (as discussed in Section~\ref{subsec:multiperiod}), most multi-period formulations in the reviewed literature capture only limited aspects of the system's dynamics, such as phased facility opening or time-varying demand. As a result, time-dependent strategic modeling remains a promising and still underexplored research direction. This limitation is particularly relevant in developing countries, where financial, structural, and political conditions are often volatile, and multiple system components evolve simultaneously.

For instance, in the Gambella region of Ethiopia discussed above, a key challenge lies in the gradual and uncertain availability of both financial and workforce resources. Decisions regarding the opening of new facilities and the deployment of qualified healthcare personnel are typically constrained by training capacity and fluctuating funding from NGOs and donors, rather than occurring simultaneously at the beginning of the planning horizon. Furthermore, forced migration driven by conflicts (e.g., the war in the neighboring country South Sudan) introduces substantial uncertainty and variability in population distribution and healthcare demand. Multi-period models that explicitly and jointly capture the temporal evolution of investment budgets, facility expansion, healthcare workforce availability, and demand fluctuations could therefore provide more realistic and actionable decision support. However, to the best of our knowledge, such comprehensive dynamic formulations remain largely absent from the existing literature.

\subsubsection{Disruptions to road links}
Due to the instability and precariousness of the healthcare systems in some developing countries, there have been a few attempts to consider disruptions to facilities (\citealp{mohammadi2014design}; \citealp{shishebori2015robust}; \citeauthor{zarrinpoor2018design},  \citeyear{zarrinpoor2017design},  \citeyear{zarrinpoor2018design}). However, only \cite{shishebori2015robust} account for disruptions to road links as well. The prevalence of adverse weather conditions and the deterioration of transportation infrastructure, increasingly exacerbated by climate change, can render certain roads impassable, preventing entire communities from reaching a particular facility, even if it remains operational \citep{van2022optimizing}. In addition, some disruptions could mean that a road normally accessible by motorized vehicles can only be traversed on foot or by boat. Modeling multi-modal modes of transport in the context of possible disruptions to road links remains an unexplored area.

\subsubsection{Patient behavior and security aspects}
Furthermore, several critical issues in developing countries have received little to no attention in the analyzed papers. For instance, factors such as hospital competition \citep{wang2018healthcare} or the reputation of different healthcare centers \citep{motallebi2020short-term} could significantly influence people's willingness to seek care at a particular facility in large developing cities.  Additionally, security aspects are of utmost importance in specific areas. For example, certain health centers might be inaccessible to particular inhabitants due to ethnic conflicts (e.g., between the Anuak and Nuar groups in the Gambella region), forcing them to attend alternative centers, even if they are farther away. Similarly, ambulance drivers may face restrictions on traveling certain roads or through specific regions. Consequently, alternative routes might need to be considered, or these drivers allocated elsewhere. Future research could incorporate constraints that completely prevent people from accessing certain facilities or roads, as well as objective functions that account for these behavior patterns using mechanisms such as penalties. While patients' preferences have been addressed in existing literature mainly in terms of travel distances (e.g., \citealp{gunecs2014}), more complex behavioral and informational barriers (such as low education levels, health literacy, or lack of trust in health services) are also prevalent and remain largely underexplored. Among the reviewed papers, only \cite{alban2022resource} explicitly incorporate such behavioral considerations by modeling effective adopter demand as a sum of sigmoidal functions that capture informational barriers and the willingness of individuals to access the services offered by mobile healthcare units. Modeling these types of concerns represents a promising step toward bringing the proposed approaches closer to reality \citep{Kunc2020}. Extending these ideas by integrating behavioral dynamics into classical location or routing models is a valuable avenue for further research.

\subsection{New optimization methods}

With respect to optimization methods, the selected papers encompass several modeling and solution techniques, including robust optimization (e.g., \citealp{shishebori2015robust}; \citealp{zarrinpoor2018design}; \citealp{motallebi2020short-term}), stochastic programming (e.g., \citealp{khodaparasti2018multi-period}; \citealp{taymaz2020healthcare}), or multi-level approaches \citep{wang2018healthcare}, among others, sometimes presenting sophisticated models that combine some of these features. Nevertheless, based on our analysis, the integration of these optimization techniques with other methodologies, such as multi-criteria decision-making (MCDM) or simulation, is scarce and has so far been addressed only tangentially, as detailed below. This observation is consistent with findings from recent reviews in the OR literature. For instance, among the healthcare-related studies in developing countries surveyed by \cite{orhan2022analytics}, only one proposes an MCDM framework and only five adopt simulation-based methodologies, none hybridized with optimization. From our perspective, this is particularly limiting in strategic healthcare planning in developing countries, where decisions typically involve multiple conflicting objectives, heterogeneous stakeholder preferences, and highly uncertain system's dynamics that cannot be fully captured with pure optimization techniques. Altogether, this points to a clear gap in the literature regarding the integration of optimization with MCDM and simulation methodologies.

With respect to MCDM, some of the analyzed works consider the preferences of decision-makers (DMs) by providing them with Pareto optimal solutions from which to choose after solving the optimization problem (e.g., \citealp{beheshtifar_multiobjective_2015}; \citealp{mousazadeh_accessible_2018}; \citealp{rouhani2022robust}), while others incorporate DMs' aspirations within the objective function using a goal programming approach (\citealp{khodaparasti2017enhancing}; \citealp{dogan_model_2020}). In highly volatile environments such as developing countries, it may be appropriate to combine MCDM interactive methods \citep{hwang2012multiple} with optimization problems, as the prospects of DMs may vary over time. Among the papers that propose multi-period optimization models, none of them explicitly integrates MCDM techniques within such dynamic frameworks. Existing multi-period models typically reflect DMs' preferences only implicitly, for example through parameter tuning or risk levels (e.g., \citealp{khodaparasti2018multi-period}), rather than through formal MCDM approaches, presenting ample opportunities for research. Furthermore, the combination of dynamic MCDM \citep{campanella2011framework} and optimization remains unexplored in these settings. As a complementary line of research, recent work illustrates how MCDM can be used to identify candidate sites for vaccine facilities, which are then evaluated within an optimization model to select optimal locations, an important but often neglected step in the facility location literature \citep{kheybari2026}. Together, these considerations highlight the need for more explicit integrations of MCDM and optimization, which would bring facility location models closer to reality and facilitate their implementation in practice.

On the other hand, hybrid optimization and simulation frameworks are scarce within the context of developing countries. To the best of our knowledge, only \cite{boutilier2020ambulance} and \cite{chouksey2022optimisation} integrate simulation in their optimization models. The former attempts to capture the effects of congestion on ambulance response times, while the latter performs a Monte Carlo simulation to analyze service levels under demand stochasticity. 
These two papers only consider variations in demand or travel times in their simulation frameworks, leaving room for numerous avenues of future research. A possibility could be to address other inherently variable features such as population growth or workforce availability from a simulation point of view within the optimization problems. In addition,  hybrid methodologies have recently been applied in healthcare systems of developed countries, showcasing improvements in both short and long-term strategic planning. In particular, \cite{Ordu2021} propose a discrete event simulation model for the UK health system that captures the entire outpatient services pathway, which later serves as an input of their optimization approach. These authors also forecast the demand for all the different specialties. \cite{Mitropoulos2023} combine system dynamics modeling with optimization to redesign the primary health system in Southern Greece, considering variable health demographics. Tailoring this type of hybrid methodologies to the specific context of developing countries would be beneficial to improve healthcare accessibility.

Regarding solution strategies, we observe a tendency to either use commercial solvers for optimally solving small- to medium-size instances or to apply evolutionary-based metaheuristics. More than half of the selected location-type papers examine case studies with no more than 100 demand points or candidate sites, with only seven works addressing instances exceeding the 1000-point threshold. However, real-world scenarios, such as optimizing healthcare systems across entire regions or countries, often involve significantly larger instances. 
For example, \cite{mendoza-gomez2022regionalization} highlight that, when considering the whole State of Mexico, a problem instance would include approximately 32,000 demand points and 3,000 potential sites. 
Given the complexity and scale of such problems, future research could explore the development of advanced methods capable of efficiently solving large-scale scenarios. In recent years, growing attention has been given to the use of machine learning (ML) techniques for solving combinatorial optimization problems, both in combination with exact methods  \citep{bengio2021} and metaheuristics \citep{karimi2022}. Within the latter category, \cite{bengio2020} proposes an ML-based method capable of providing fast and good quality initial solutions to a facility location problem with uncertain demand. Based on the findings of this review, such approaches could be particularly valuable for application in developing countries. 

In the context of multi-objective problems, most selected papers rely on the classical weighted sum or $\varepsilon$-constraint methods, as well as heuristic approaches. Multi-objective optimization has recently evolved with the proposal of more sophisticated techniques, such as developing branch-and-bound methods \citep{Adelgren2022}, and providing the DM with a representative set of the Pareto front (\citealp{Dogan2022}; \citealp{MesquitaCunha2023}). Given the conflicting objectives in improving access to healthcare (especially when comparing patient and system perspectives), these novel approaches should be explored in future studies. Another promising direction is the use of enhanced variants of evolutionary-based metaheuristics. Among the selected papers, the non-dominated sorting genetic algorithm II (NSGA-II, \citealp{deb2002}) was applied in three studies (\citealp{beheshtifar_multiobjective_2015}; \citealp{motallebi2020short-term}; \citealp{decampos2024multi}), while its enhanced version, NSGA-III, \citep{deb2013}, and more recent improvements \citep{Gu2020} remain unexplored.

\subsection{Tools for ensuring practical implementation}

Mirroring the findings of studies in other application areas such as humanitarian logistics \citep{DeVries2020}, our analysis highlights that the practical adoption of optimization models to improve health service accessibility in developing countries is still scant, with only four studies documenting some level of real implementation. It also reveals some of the main challenges to practical implementation, which include: limited data availability, difficulty of accurate parameters' estimation, plurality of stakeholders with different expectations and goals, lack of local knowledge, communication barriers, and skepticism to model adoption. Many of these challenges are interlinked and attributable to the lack of effective engagement with key stakeholders throughout the model development. Furthermore, challenges such as data limitations or stakeholder plurality are common across OR applications, but their severity is often amplified in developing countries due to weak infrastructure, unstable governance, fragmented funding mechanisms, and limited technical capacity at the local level, making practical implementation even more difficult.

Future research should aim at redressing this shortcoming by developing multi-methodology approaches \citep{Mingers1997} where optimization is combined with problem structuring methods (PSMs) to support health decision-making in practice (see for example \citealp{Cardoso2019}). Importantly, the role of PSMs in this context is not only to foster acceptance of OR models, but also to help determine whether, where, and how optimization can meaningfully add value given the practical constraints of the setting under study. PSMs, such as Soft Systems Methodology (SSM), Strategic Options Development and Analysis (SODA), and Strategic Choice Approach (SCA), are a set of rigorous qualitative methods specifically designed to address complex, ``wicked'' problems with multiple stakeholders and conflicting interests \citep{PetropoulosEtAl2024}. These methods focus on framing issues, facilitating group discussions, and developing a shared understanding of problematic situations, with a view to reaching a consensus on problem structure and fostering stakeholder agreement to act (see \citealp{SmithShaw2019, Franco2022} for recent reviews). 

In the complex healthcare settings of developing countries, stakeholders include a wide range of actors such as government health officials, local healthcare providers, NGOs, community leaders, and patients themselves, each bringing unique perspectives and priorities. PSMs offer structured frameworks that facilitate inclusive dialogue among these diverse groups. Through this process, OR practitioners can gain insight into the realistic nature of the practical strategic planning problems they are facing.  
These methods help capture the conflicting objectives and constraints of stakeholders, such as limited resources, cultural barriers, and logistical challenges, enabling the development of optimization models that are grounded in the practical realities of the local context.
The adoption of these participatory approaches would ensure that models reflect the needs and priorities of those directly impacted by healthcare accessibility issues, such as rural populations, marginalized communities, and healthcare workers. They would also help translate technical outputs into actionable strategies by incorporating local knowledge and ensuring that solutions are both culturally sensitive and practically feasible. In addition, they also help avoid over-modeling situations where the feasible decision space is too limited for optimization to provide significant benefits over simpler decision-support approaches. 

Despite their enormous potential to address grand challenges \citep{Abuabara2021, Ackermann2024} and the advocacy that PSMs should be an essential element of any OR engagement \citep{Dyson2021}, PSMs have not yet been proposed to tackle healthcare accessibility problems in developing countries, offering new research opportunities. Within this specific context, a way forward is the adoption of PartiOpt, a participative framework recently introduced in \cite{noorain2024} to actively engage stakeholders in the development of healthcare optimization models \citep{Noorain2025}. PartiOpt builds on PartiSim \citep{tako2015}, a widely used multi-methodology framework for discrete event simulation, and offers tools for stakeholder engagement across all stages of optimization model development, from problem definition to model validation, experimentation, and implementation. 


PartiOpt also includes structured tools to identify data needs, availability, accessibility, and ownership. As evidenced in the four papers involving practical implementation, the use of real data is key to ensuring model adoption. Gathering health system data in developing countries is a challenging task and often requires using a range of tools, as demonstrated in the real case study by \cite{de2020roadside}. As seen in the analyzed papers, geospatial data are typically obtained through GIS software and GPS devices. Global databases also exist that store a variety of health data, such as disease burden, service availability, health program indicators, vaccination, mortality, and morbidity statistics. Noteworthy examples include the WHO Global Health Observatory (\url{https://www.who.int/data/gho}) and the District Health Information System 2 (\url{https://dhis2.org/}). However, despite substantial progress in global health data collection, developing countries still face significant data gaps and challenges. Available health data are often aggregated at the district or national level without sufficient disaggregation to reveal local health disparities and needs, especially in rural, marginalized, or conflict-affected areas. In addition, critical data on health infrastructure, workforce, and patient condition or behavior are often missing.  
As a result, the practical implementation of optimization models may require tailored and context-specific primary data collection using multiple tools, ranging from mobile technology and GPS for real-time patient flow and service utilization, to household surveys, health expert interviews, and focus group workshops to elicit views on parameter estimates, candidate facility locations, community priorities, and behavioral aspects. Organizing and financing this primary data collection generally requires the collaboration among local health authorities, NGOs, and international organizations, such as the UN. This phase is critical for testing, validating, and implementing optimization models. Without reliable real-world data, the applicability and impact of models are naturally limited, and recommendations derived from theoretical or synthetic data may be misleading. 

To enhance the realism and applicability of optimization models, artificial intelligence (AI) tools could also be used for facilitating health data integration and analysis, as well as to generate inputs for the models, such as health trends and disease outbreak predictions. The use of AI, including ML and deep learning, would also be beneficial to guide health system decision-making in regions prone to frequent natural disasters. As an example, AI can be used to generate hazard zoning maps for areas affected by floods to classify zones based on their risk levels \citep{Singha2024}. These maps can then inform decisions on new health facilities' locations and transport infrastructure improvements.               

Finally, as with other sustainable development challenges that involve complex social, economic, geographic, and cultural factors, developing effective optimization models for healthcare accessibility may require an interdisciplinary approach combining data, methods, and findings from various fields (see \citealp{scaparra2024}, for an example in disaster management). Relevant disciplines include public health, economics, and geography to provide inputs about disease incidence, patient demographics, health delivery models, cost structures, and existing infrastructure. Sociology and anthropology can help understand how cultural values and human factors, such as trust in service providers, influence health-seeking behaviors and service utilization. Additionally, social marketing can influence the demand side of healthcare access by encouraging behavioral changes and raising awareness about service provision through culturally sensitive communication strategies. This synergy with optimization models would help ensure not only an improved supply of services but also greater acceptance and use.

\section{Conclusions}\label{sec:conclusions}

This review classifies and analyzes optimization-based approaches proposed over the past two decades to enhance geographical accessibility to health services in developing countries. The surveyed works are categorized by problem characteristics, optimization methodologies, and practical implementation. A roadmap is also provided to highlight research gaps and propose future research directions.

In terms of problem features, most studies concentrate on locating primary care fixed facilities and aim to maximize coverage, minimize patient travel, or reduce system costs. Accessibility, equity, and efficiency are the predominant high-level criteria, and only about 29\%,  34\%, and 21\% of the papers incorporate models that account for the hierarchical structure, uncertainty, or multi-period aspects of health systems, respectively. Furthermore, case studies are unevenly distributed, with Iran and India overrepresented, likely due to their strong OR communities. Future research should explore specialized services, workforce allocation, mobile and fixed facility integration, comprehensive multi-period settings, and factors like road disruptions, patient behavior, and security. Expanding research efforts in underrepresented regions, particularly Sub-Saharan Africa and South-East Asia, is also crucial for achieving broader progress in healthcare accessibility.

Methodologically, most papers address location-allocation problems modeled as single-objective integer or mixed-integer linear programs. A range of methods like fuzzy programming, stochastic programming, robust optimization, and goal programming are also employed. Common techniques to handle multiple objectives are the weighted sum method, the $\varepsilon$-constraint method, and evolutionary algorithms. In terms of solution approaches, exact methods implemented via commercial solvers dominate, with only a few studies using ad-hoc exact methods or stand-alone approximate methods. We argue that integrating optimization with other techniques, such as MCDM, simulation, or ML models, could significantly enrich the OR literature in this field.

The practical implementation section examines the case studies, focusing on data sources, stakeholder engagement, and real-world application. Most models remain theoretical proofs of concept, with only four studies reporting actual implementation, likely due to the limited involvement of real stakeholders in the design and progression of the research projects. Integrating optimization with PSMs could improve engagement and applicability. Additionally, multidisciplinary approaches are crucial for enhancing model realism and directly influencing the health-seeking behavior patterns of the population. Furthermore, without access to real-world data, the quality of the results obtained and the reliability of the recommendations from the models are uncertain, which underscores the importance of structured data collection as an integral part of the research process.

Finally, achieving lasting improvements in healthcare accessibility requires investing in capacity building and collaboration across different roles. The development of mathematical models and solution methods should primarily be carried out by OR experts, in close collaboration with computer scientists or software engineers who can translate these models into user-friendly decision-support tools. NGO staff and local government officials should then be trained to effectively use, interpret, and trust these tools in their decision-making processes. This would definitely enhance the reliability and usability of optimization-based solutions, and facilitate their effective adoption in practice. The focus on academic publications often overlooks the complexities of real-world adoption. Fostering sustained collaboration between academia, software developers, and stakeholders is therefore essential to bridge the gap between theory and real impact, ultimately improving healthcare outcomes in developing countries.

\section*{Acknowledgments}
L. Davila-Pena acknowledges the support of the R+D+I project PID2021-124030NB-C32, granted by MICIU/AEI/10.13039/501100011033/ and by ``ERDF A way of making Europe''/EU. She also acknowledges support from the Xunta de Galicia (Grupos de Referencia Competitiva ED431C 2025/03). The authors are grateful to Prof. Joaquim Gromicho, Prof. Goos Kant, and Ms. Parvathy Krishnan Krishnakumari for their constructive feedback in the early stages of this research. 
 Finally, a very special thanks to Doctors with Africa CUAMM for their invaluable insights.




\bibliography{references_black}

\include{supplementary_material_black}

\end{document}

%% file: supplementary_material_black.tex
%
%
%
%
%
%
%
%
%
%
%

\clearpage

\setcounter{page}{1}
\setcounter{section}{0}
\setcounter{table}{0}
\setcounter{figure}{0}
\setcounter{equation}{0}

\renewcommand{\thepage}{S\arabic{page}}
\renewcommand{\thesection}{S\arabic{section}}
\renewcommand{\thetable}{S\arabic{table}}
\renewcommand{\thefigure}{S\arabic{figure}}
\renewcommand{\theequation}{S\arabic{equation}}

\begin{center}
	 \Large{\textbf{Supplementary Material}}
\end{center}

\section{Additional details on search strategy}\label{sec:strategy}

The literature search was conducted across INFORMS journals, ScienceDirect, and Google Scholar databases, following a three-step strategy. In the first step, searches were performed in the INFORMS journals and ScienceDirect databases. For the INFORMS journals database, an advanced search was conducted using specific Boolean search strings to retrieve works that align with the review's objectives. These strings targeted keywords related to accessibility, Operations Research, healthcare, and developing countries within the designated timeframe (January 2003--June 2025), with searches spanning across the entire documents. Only research articles written in English were considered for inclusion. 
Due to limitations in the number of Boolean operators on ScienceDirect, a refined search strategy was implemented. This involved variations in search terms to ensure comprehensive coverage of relevant literature. Table~\ref{tab:databases} provides details on specific terms and parameters that were used. Adjustments in ScienceDirect included substituting ``\textit{operations research}'' with ``\textit{optimization}'', ``\textit{access}'' with ``\textit{barrier}'' or ``\textit{coverage}'', incorporating alternative terms for ``\textit{developing countries}'' such as ``\textit{developing nations}'', removing/including hyphens, or considering UK English spelling variations, among others.

\begin{table}[H]
    \centering \footnotesize
    \caption{Databases and search strategy.}
    \label{tab:databases}
    \begin{tabular}{c p{11cm}}
        \noalign{\smallskip}\hline\noalign{\smallskip}
        \textbf{Database} & \textbf{Search criteria} \\
        \noalign{\smallskip}\hline\noalign{\smallskip}
         & Anywhere: (“\textit{accessibility}” OR “\textit{access}”)\\ 
          & AND\\ 
        & Anywhere: (“\textit{operations research}” OR “\textit{operational research}” OR “\textit{optimization}” OR “\textit{optimisation}”)\\
        INFORMS & AND\\
        journals & Anywhere: (“\textit{developing countries}” OR “\textit{low-income countries}” OR “\textit{low- and middle-income countries}” OR “\textit{low and middle income countries}” OR  “\textit{LMIC}” OR “\textit{low-middle income countries}” OR “\textit{lesser developed countries}” OR “\textit{lesser-developed countries}”)\\
        & AND \\
        & Anywhere: (“\textit{healthcare}” OR “\textit{health care}” OR “\textit{health}”) \\
        \noalign{\smallskip}\hline\noalign{\smallskip}
        \multirow{6}{*}{ScienceDirect} & Articles with these terms: (“\textit{accessibility}” OR “\textit{access}”) AND (“\textit{operations research}”) AND (“\textit{developing countries}” OR “\textit{low income countries}” OR “\textit{LMIC}” OR “\textit{low and middle income countries}” OR “\textit{lesser developed countries}” OR “\textit{underdeveloped countries}”)\\
        & AND\\
        & Title, abstract, keywords: (“\textit{healthcare}” OR “\textit{health care}” OR “\textit{health}”) \\
        \noalign{\smallskip}\hline\noalign{\smallskip}
    \end{tabular}
\end{table}

This first step resulted in the inclusion of 13 papers that met the review's eligibility criteria. In the second step, these papers were searched in the Scopus database, and backward and forward citation searches were conducted for each of them. This process led to the identification of 30 additional papers. In the third step, a manual backward search was performed on these newly identified papers, resulting in the inclusion of 13 further records. Overall, a total of 56 papers were included in the final review. 

\section{Brief bibliometric analysis}\label{sec:bibliometric_analysis}
\begin{figure}[ht!]
    \centering
    \includegraphics[width=\linewidth]{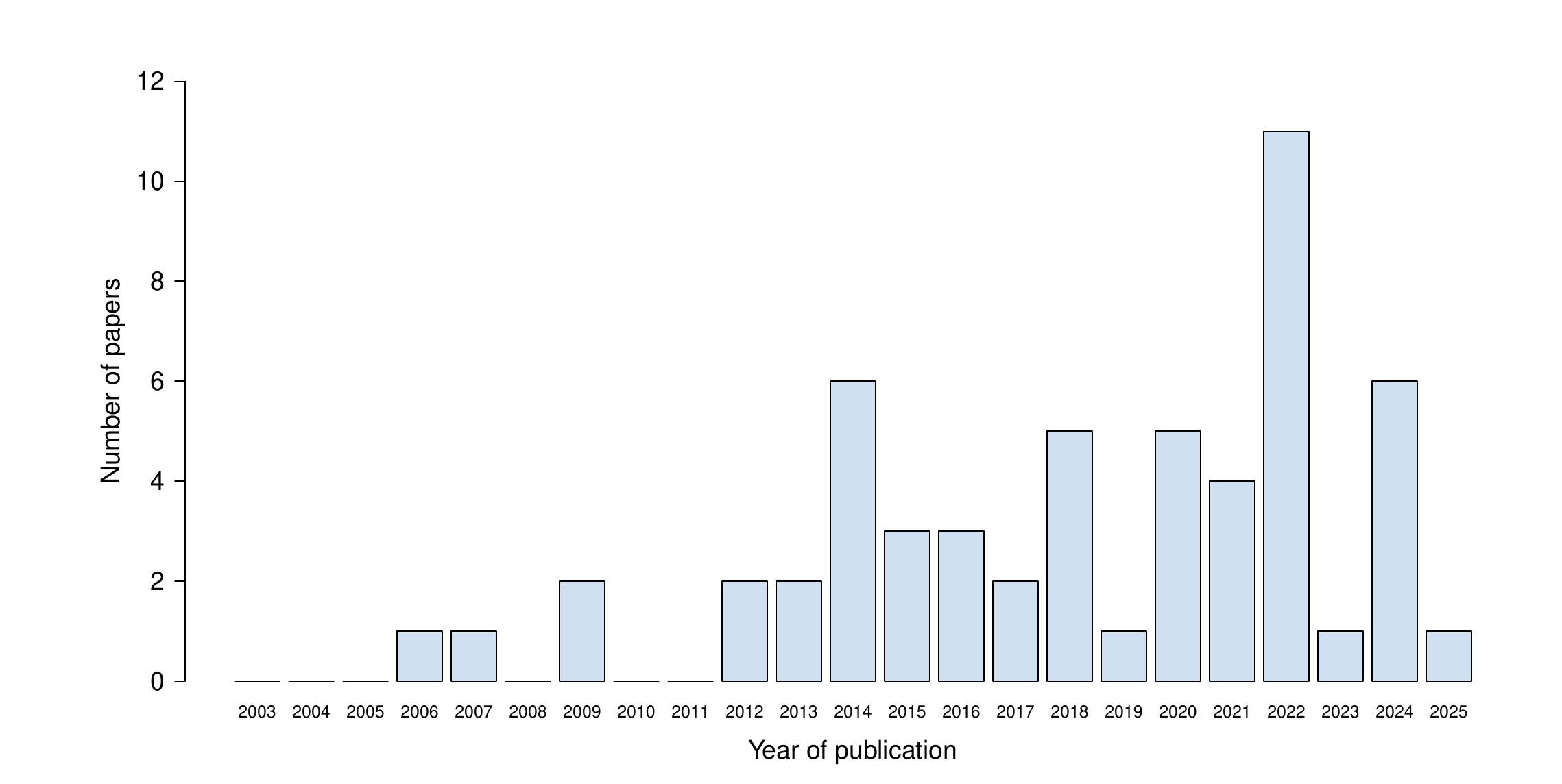}
    \caption{Number of selected papers by year of publication.}
    \label{fig:barplot_year}
\end{figure}

\begin{figure}[ht!]
    \centering
    \includegraphics[width=\linewidth]{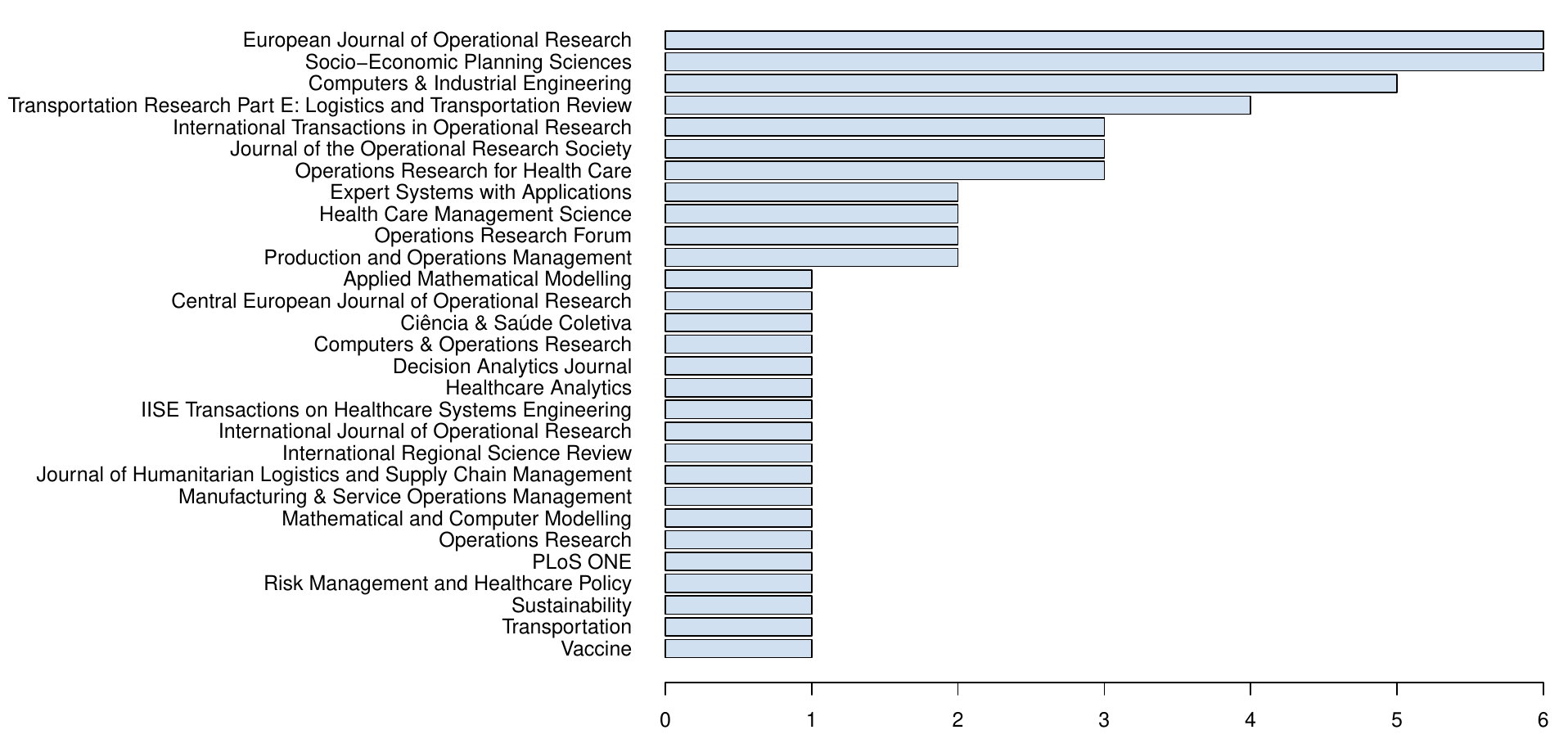} 
    \caption{Number of selected papers by journal.}
    \label{fig:barplot_journal}
\end{figure}

\begin{figure}[ht!]
    \centering
    \includegraphics[width=\linewidth]{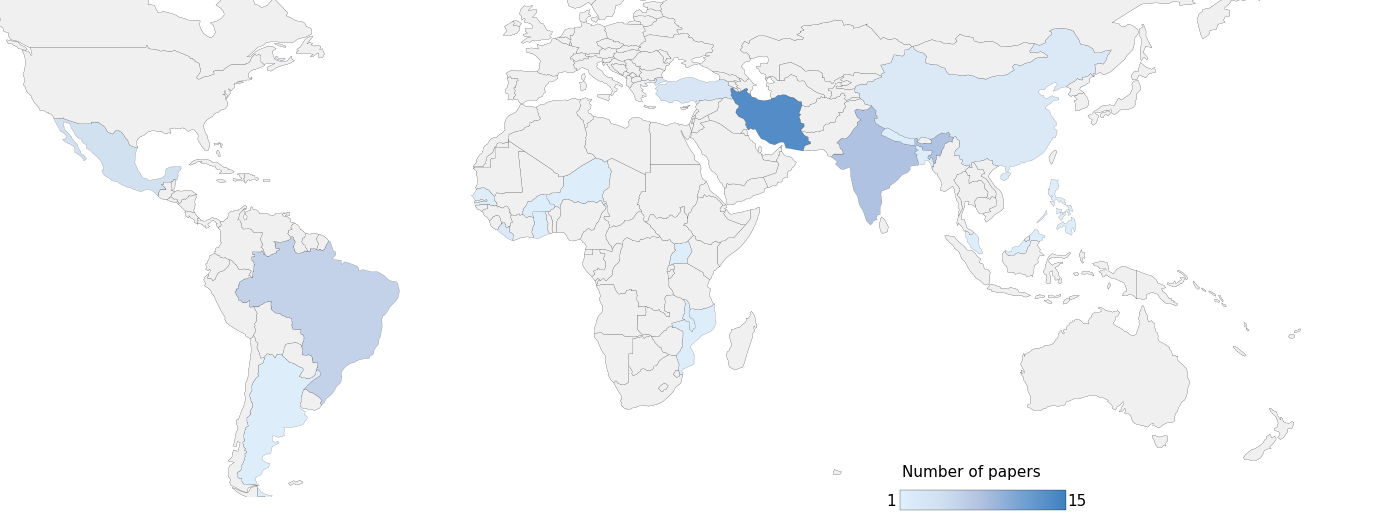}
    \caption{Countries with case study.}
    \label{fig:map}
\end{figure}

Figures~\ref{fig:barplot_year} and \ref{fig:barplot_journal} show the number of selected papers published by year and by journal. Also, Figure~\ref{fig:map} illustrates the number of publications by country where the case study takes place. Note that the map does not display those four papers with case studies in unspecified countries of Sub-Saharan and South-East Africa. 

From Figure~\ref{fig:barplot_year}, we observe that none of the selected papers was published in the periods 2003-2005 and 2010-2011, with the majority emerging within the last half of the timeframe (over the 80\% of the total papers). This trend suggests a growing interest in this topic.

Additionally, Figure~\ref{fig:barplot_journal} shows that the European Journal of Operational Research (EJOR) and Socio-Economic Planning Sciences are the two leading journals in terms of publication volume. Notably, most of the papers are published in OR-focused journals, with a few exceptions found in specialized fields such as geography or medicine.

As we can see on the map in Figure~\ref{fig:map}, Iran accounts for over 26\% of the case studies, followed by India with approximately 14\%. This is not surprising, given the presence of strong Operations Research (OR) communities in these countries. We also observe that among all the developing countries in the Americas, only the largest ones (Brazil, Argentina, and Mexico) have been studied six, one, or four times, respectively. While some studies focus on low- and middle-income countries in Africa and Asia, these continents remain largely unexplored.

We have also examined whether some authors had affiliations in the countries where the case studies were conducted, revealing 17 studies without such affiliations and 39 with them. Excluding cases from Iran and India, only 15 of the remaining 33 papers feature authors affiliated with institutions in the respective case study countries. Notably, none of the 13 papers containing case studies in African countries are authored by individuals affiliated with that continent.
